\documentclass[12pt,a4paper]{amsart}
\usepackage[utf8]{inputenc}%
\usepackage[T1]{fontenc}%
\usepackage[english]{babel}%
\usepackage[margin=2cm]{geometry}%
\usepackage{dsfont}%
\usepackage[pdftex]{graphicx}%
\usepackage{amsmath,amssymb,amsthm}%
\usepackage{color}%
\usepackage{float}%
\usepackage{ifthen}%
\usepackage{hyperref}%
\usepackage{tikz}
\usepackage{pgfplots}
\usepackage{mathtools}
\usepackage{booktabs}
\usepackage{pdfpages}
\usepackage[onehalfspacing]{setspace}
\usepackage{amssymb,amsmath,amsfonts}
\usepackage{amsthm}
\usepackage{doi}
\usepackage{algorithm}
\usepackage{algpseudocode}
\algrenewcommand\algorithmicrequire{\textbf{Input:}}
\algrenewcommand\algorithmicensure{\textbf{Output:}}
\usepackage{hyperref}
\hypersetup{
    hidelinks,
    plainpages=false
}
\usepackage{microtype}
\usepackage{tcolorbox}
\usepackage{enumitem}
\usepackage{longtable}
\usepackage{csquotes}%
\usepackage[labelformat=simple]{subcaption}

\pgfplotsset{compat=newest}
\usepgfplotslibrary{external} 
\usepackage{mathtools}
\usepgfplotslibrary{colormaps}
\usepgfplotslibrary{groupplots} %
\usetikzlibrary{pgfplots.groupplots} %
\usepackage{multirow}
\usetikzlibrary{decorations.markings}
\pgfplotsset{compat=1.18}

\newcommand*{\rom}[1]{\expandafter\@slowromancap\romannumeral #1@}

\definecolor{carrotorange}{rgb}{0.93, 0.57, 0.13}%
\definecolor{kellygreen}{rgb}{0.3, 0.73, 0.09}%
\definecolor{brightcerulean}{rgb}{0.11, 0.67, 0.84} %
\definecolor{cardinal}{rgb}{0.77, 0.12, 0.23} %
\definecolor{tablegray}{gray}{0.95}
\definecolor{dpgreen}{rgb}{0.01, 0.76, 0.23}%
\definecolor{forestgreen}{rgb}{0.13, 0.55, 0.13}

\def\for{\mathrm{for~}}

\def\and{\mathrm{and~}}

\newcommand{\VSBDF}[1]{\textnormal{VS-BDF}\ensuremath{(#1)}}
\newcommand{\BDF}[1]{\textnormal{BDF}\ensuremath{(#1)}}

\newcommand{\noopsort}[1]{}

\def\qq{\qquad}
\def\qqq{\qquad\qq}

\newcommand*\dss{\mathop{}\!\mathrm{d}s}
\newcommand*\dt{\mathop{}\!\mathrm{d}t}

\newcommand{\R}{\mathbb{R}} %
\renewcommand{\P}{\mathbb{P}} %
\newcommand{\N}{\mathbb{N}} %

\newcommand{\norm}[1]{\lVert #1 \rVert}

\newcommand{\dual}[1]{\langle #1 \rangle}

\def\grad{\boldsymbol{\nabla}}%
\newcommand{\jump}[1]{ [\![#1]\!]_E}

\newcommand{\recM}[1]{U_h^{(#1)}}
\newcommand{\recW}[1]{W^{(#1)}}
\newcommand{\errM}[1]{\varepsilon^{(#1)}}
\newcommand{\errW}[1]{\tilde{\varepsilon}^{(#1)}}
\newcommand{\recMn}[2]{U_h^{{#1,}(#2)}}
\newcommand{\recWn}[2]{W^{{#1,}(#2)}}

\def\ee{\mathrm{e}}%
\def\ptl{\partial}%

\renewcommand{\bar}{\overline}

\renewcommand{\tilde}{\widetilde}

\theoremstyle{plain}%
\newtheorem{Theorem}{Theorem}

\newtheorem{Lemma}[Theorem]{Lemma}
\newtheorem{Corollary}[Theorem]{Corollary}
\newtheorem{Definition}[Theorem]{Definition}
\newtheorem{Remark}[Theorem]{Remark}

\newtheorem{Assumption}[Theorem]{Assumption}
\numberwithin{equation}{section}
\numberwithin{table}{section}
\numberwithin{figure}{section}

\def\EE{{\mathcal E}}

\def\II{{\mathcal I}}

\def\NN{{\mathcal N}}

\def\RR{{\mathcal R}}
\def\SS{{\mathcal S}}

\def\TT{{\mathcal T}}
\def\FF{{\mathcal F}}

\def\bL{\boldsymbol{L}}%
\def\bn{\boldsymbol{n}}%
\def\b0{\boldsymbol{0}}%

\definecolor{kit-green}{RGB}{0, 150, 130}
\definecolor{kit-green100}{RGB}{0, 150, 130}
\definecolor{kit-green90}{rgb}{0.1, 0.6294, 0.5588}
\definecolor{kit-green80}{rgb}{0.2, 0.6706, 0.6078}
\definecolor{kit-green75}{rgb}{0.25, 0.6912, 0.6324}
\definecolor{kit-green70}{rgb}{0.3, 0.7118, 0.6569}
\definecolor{kit-green60}{rgb}{0.4, 0.7529, 0.7059}
\definecolor{kit-green50}{rgb}{0.5, 0.7941, 0.7549}
\definecolor{kit-green40}{rgb}{0.6, 0.8353, 0.8039}
\definecolor{kit-green30}{rgb}{0.7, 0.8765, 0.8529}
\definecolor{kit-green25}{rgb}{0.75, 0.8971, 0.8775}
\definecolor{kit-green20}{rgb}{0.8, 0.9176, 0.902}
\definecolor{kit-green15}{rgb}{0.85, 0.9382, 0.9265}
\definecolor{kit-green10}{rgb}{0.9, 0.9588, 0.951}
\definecolor{kit-green5}{rgb}{0.95, 0.9794, 0.9755}

\definecolor{kit-blue}{RGB}{70, 100, 170}
\definecolor{kit-blue100}{RGB}{70, 100, 170}
\definecolor{kit-blue90}{rgb}{0.3471, 0.4529, 0.7}
\definecolor{kit-blue80}{rgb}{0.4196, 0.5137, 0.7333}
\definecolor{kit-blue75}{rgb}{0.4559, 0.5441, 0.75}
\definecolor{kit-blue70}{rgb}{0.4922, 0.5745, 0.7667}
\definecolor{kit-blue60}{rgb}{0.5647, 0.6353, 0.8}
\definecolor{kit-blue50}{rgb}{0.6373, 0.6961, 0.8333}
\definecolor{kit-blue40}{rgb}{0.7098, 0.7569, 0.8667}
\definecolor{kit-blue30}{rgb}{0.7824, 0.8176, 0.9}
\definecolor{kit-blue25}{rgb}{0.8186, 0.848, 0.9167}
\definecolor{kit-blue20}{rgb}{0.8549, 0.8784, 0.9333}
\definecolor{kit-blue15}{rgb}{0.8912, 0.9088, 0.95}
\definecolor{kit-blue10}{rgb}{0.9275, 0.9392, 0.9667}
\definecolor{kit-blue5}{rgb}{0.9637, 0.9696, 0.9833}

\definecolor{kit-red}{RGB}{162, 34, 35}
\definecolor{kit-red100}{RGB}{162, 34, 35}
\definecolor{kit-red90}{rgb}{0.6718, 0.22, 0.2235}
\definecolor{kit-red80}{rgb}{0.7082, 0.3067, 0.3098}
\definecolor{kit-red75}{rgb}{0.7265, 0.35, 0.3529}
\definecolor{kit-red70}{rgb}{0.7447, 0.3933, 0.3961}
\definecolor{kit-red60}{rgb}{0.7812, 0.48, 0.4824}
\definecolor{kit-red50}{rgb}{0.8176, 0.5667, 0.5686}
\definecolor{kit-red40}{rgb}{0.8541, 0.6533, 0.6549}
\definecolor{kit-red30}{rgb}{0.8906, 0.74, 0.7412}
\definecolor{kit-red25}{rgb}{0.9088, 0.7833, 0.7843}
\definecolor{kit-red20}{rgb}{0.9271, 0.8267, 0.8275}
\definecolor{kit-red15}{rgb}{0.9453, 0.87, 0.8706}
\definecolor{kit-red10}{rgb}{0.9635, 0.9133, 0.9137}
\definecolor{kit-red5}{rgb}{0.9818, 0.9567, 0.9569}

\definecolor{kit-yellow}{RGB}{252, 229, 0}
\definecolor{kit-yellow100}{RGB}{252, 229, 0}
\definecolor{kit-yellow90}{rgb}{0.9894, 0.9082, 0.1}
\definecolor{kit-yellow80}{rgb}{0.9906, 0.9184, 0.2}
\definecolor{kit-yellow75}{rgb}{0.9912, 0.9235, 0.25}
\definecolor{kit-yellow70}{rgb}{0.9918, 0.9286, 0.3}
\definecolor{kit-yellow60}{rgb}{0.9929, 0.9388, 0.4}
\definecolor{kit-yellow50}{rgb}{0.9941, 0.949, 0.5}
\definecolor{kit-yellow40}{rgb}{0.9953, 0.9592, 0.6}
\definecolor{kit-yellow30}{rgb}{0.9965, 0.9694, 0.7}
\definecolor{kit-yellow25}{rgb}{0.9971, 0.9745, 0.75}
\definecolor{kit-yellow20}{rgb}{0.9976, 0.9796, 0.8}
\definecolor{kit-yellow15}{rgb}{0.9982, 0.9847, 0.85}
\definecolor{kit-yellow10}{rgb}{0.9988, 0.9898, 0.9}
\definecolor{kit-yellow5}{rgb}{0.9994, 0.9949, 0.95}

\definecolor{kit-orange}{RGB}{223, 155, 27}
\definecolor{kit-orange100}{RGB}{223, 155, 27}
\definecolor{kit-orange90}{rgb}{0.8871, 0.6471, 0.1953}
\definecolor{kit-orange80}{rgb}{0.8996, 0.6863, 0.2847}
\definecolor{kit-orange75}{rgb}{0.9059, 0.7059, 0.3294}
\definecolor{kit-orange70}{rgb}{0.9122, 0.7255, 0.3741}
\definecolor{kit-orange60}{rgb}{0.9247, 0.7647, 0.4635}
\definecolor{kit-orange50}{rgb}{0.9373, 0.8039, 0.5529}
\definecolor{kit-orange40}{rgb}{0.9498, 0.8431, 0.6424}
\definecolor{kit-orange30}{rgb}{0.9624, 0.8824, 0.7318}
\definecolor{kit-orange25}{rgb}{0.9686, 0.902, 0.7765}
\definecolor{kit-orange20}{rgb}{0.9749, 0.9216, 0.8212}
\definecolor{kit-orange15}{rgb}{0.9812, 0.9412, 0.8659}
\definecolor{kit-orange10}{rgb}{0.9875, 0.9608, 0.9106}
\definecolor{kit-orange5}{rgb}{0.9937, 0.9804, 0.9553}

\definecolor{kit-lightgreen}{RGB}{140, 182, 60}
\definecolor{kit-lightgreen100}{RGB}{140, 182, 60}
\definecolor{kit-lightgreen90}{rgb}{0.5941, 0.7424, 0.3118}
\definecolor{kit-lightgreen80}{rgb}{0.6392, 0.771, 0.3882}
\definecolor{kit-lightgreen75}{rgb}{0.6618, 0.7853, 0.4265}
\definecolor{kit-lightgreen70}{rgb}{0.6843, 0.7996, 0.4647}
\definecolor{kit-lightgreen60}{rgb}{0.7294, 0.8282, 0.5412}
\definecolor{kit-lightgreen50}{rgb}{0.7745, 0.8569, 0.6176}
\definecolor{kit-lightgreen40}{rgb}{0.8196, 0.8855, 0.6941}
\definecolor{kit-lightgreen30}{rgb}{0.8647, 0.9141, 0.7706}
\definecolor{kit-lightgreen25}{rgb}{0.8873, 0.9284, 0.8088}
\definecolor{kit-lightgreen20}{rgb}{0.9098, 0.9427, 0.8471}
\definecolor{kit-lightgreen15}{rgb}{0.9324, 0.9571, 0.8853}
\definecolor{kit-lightgreen10}{rgb}{0.9549, 0.9714, 0.9235}
\definecolor{kit-lightgreen5}{rgb}{0.9775, 0.9857, 0.9618}

\definecolor{kit-purple}{RGB}{163, 16, 124}
\definecolor{kit-purple100}{RGB}{163, 16, 124}
\definecolor{kit-purple90}{rgb}{0.6753, 0.1565, 0.5376}
\definecolor{kit-purple80}{rgb}{0.7114, 0.2502, 0.589}
\definecolor{kit-purple75}{rgb}{0.7294, 0.2971, 0.6147}
\definecolor{kit-purple70}{rgb}{0.7475, 0.3439, 0.6404}
\definecolor{kit-purple60}{rgb}{0.7835, 0.4376, 0.6918}
\definecolor{kit-purple50}{rgb}{0.8196, 0.5314, 0.7431}
\definecolor{kit-purple40}{rgb}{0.8557, 0.6251, 0.7945}
\definecolor{kit-purple30}{rgb}{0.8918, 0.7188, 0.8459}
\definecolor{kit-purple25}{rgb}{0.9098, 0.7657, 0.8716}
\definecolor{kit-purple20}{rgb}{0.9278, 0.8125, 0.8973}
\definecolor{kit-purple15}{rgb}{0.9459, 0.8594, 0.9229}
\definecolor{kit-purple10}{rgb}{0.9639, 0.9063, 0.9486}
\definecolor{kit-purple5}{rgb}{0.982, 0.9531, 0.9743}

\definecolor{kit-brown}{RGB}{167, 130, 46}
\definecolor{kit-brown100}{RGB}{167, 130, 46}
\definecolor{kit-brown90}{rgb}{0.6894, 0.5588, 0.2624}
\definecolor{kit-brown80}{rgb}{0.7239, 0.6078, 0.3443}
\definecolor{kit-brown75}{rgb}{0.7412, 0.6324, 0.3853}
\definecolor{kit-brown70}{rgb}{0.7584, 0.6569, 0.4263}
\definecolor{kit-brown60}{rgb}{0.7929, 0.7059, 0.5082}
\definecolor{kit-brown50}{rgb}{0.8275, 0.7549, 0.5902}
\definecolor{kit-brown40}{rgb}{0.862, 0.8039, 0.6722}
\definecolor{kit-brown30}{rgb}{0.8965, 0.8529, 0.7541}
\definecolor{kit-brown25}{rgb}{0.9137, 0.8775, 0.7951}
\definecolor{kit-brown20}{rgb}{0.931, 0.902, 0.8361}
\definecolor{kit-brown15}{rgb}{0.9482, 0.9265, 0.8771}
\definecolor{kit-brown10}{rgb}{0.9655, 0.951, 0.918}
\definecolor{kit-brown5}{rgb}{0.9827, 0.9755, 0.959}

\definecolor{kit-cyan}{RGB}{35, 161, 224}
\definecolor{kit-cyan100}{RGB}{35, 161, 224}
\definecolor{kit-cyan90}{rgb}{0.2235, 0.6682, 0.8906}
\definecolor{kit-cyan80}{rgb}{0.3098, 0.7051, 0.9027}
\definecolor{kit-cyan75}{rgb}{0.3529, 0.7235, 0.9088}
\definecolor{kit-cyan70}{rgb}{0.3961, 0.742, 0.9149}
\definecolor{kit-cyan60}{rgb}{0.4824, 0.7788, 0.9271}
\definecolor{kit-cyan50}{rgb}{0.5686, 0.8157, 0.9392}
\definecolor{kit-cyan40}{rgb}{0.6549, 0.8525, 0.9514}
\definecolor{kit-cyan30}{rgb}{0.7412, 0.8894, 0.9635}
\definecolor{kit-cyan25}{rgb}{0.7843, 0.9078, 0.9696}
\definecolor{kit-cyan20}{rgb}{0.8275, 0.9263, 0.9757}
\definecolor{kit-cyan15}{rgb}{0.8706, 0.9447, 0.9818}
\definecolor{kit-cyan10}{rgb}{0.9137, 0.9631, 0.9878}
\definecolor{kit-cyan5}{rgb}{0.9569, 0.9816, 0.9939}

\definecolor{kit-gray}{RGB}{0, 0, 0}
\definecolor{kit-gray100}{RGB}{0, 0, 0}
\definecolor{kit-gray90}{rgb}{0.1, 0.1, 0.1}
\definecolor{kit-gray80}{rgb}{0.2, 0.2, 0.2}
\definecolor{kit-gray75}{rgb}{0.25, 0.25, 0.25}
\definecolor{kit-gray70}{rgb}{0.3, 0.3, 0.3}
\definecolor{kit-gray60}{rgb}{0.4, 0.4, 0.4}
\definecolor{kit-gray50}{rgb}{0.5, 0.5, 0.5}
\definecolor{kit-gray40}{rgb}{0.6, 0.6, 0.6}
\definecolor{kit-gray30}{rgb}{0.7, 0.7, 0.7}
\definecolor{kit-gray25}{rgb}{0.75, 0.75, 0.75}
\definecolor{kit-gray20}{rgb}{0.8, 0.8, 0.8}
\definecolor{kit-gray15}{rgb}{0.85, 0.85, 0.85}
\definecolor{kit-gray10}{rgb}{0.9, 0.9, 0.9}
\definecolor{kit-gray5}{rgb}{0.95, 0.95, 0.95}

\newcommand{\alignedintertext}[1]{%
  \noalign{%
    \vskip\belowdisplayshortskip
    \vtop{\hsize=\linewidth#1\par
    \expandafter}%
    \expandafter\prevdepth\the\prevdepth
  }%
}

\title[A posteriori error estimates for higher-order BDF methods]{A posteriori error estimates for higher-order BDF methods for parabolic partial differential equations}

\author{Stefan Karch}
\address{TU Wien, Institute of Analysis and Scientific Computing, Wiedner Hauptstr. 8-10, 1040 Wien, Austria}
\email{stefan.karch@asc.tuwien.ac.at}

\begin{document}
\clearpage

\begin{abstract}
We derive \emph{a~posteriori} error estimates for parabolic problems discretized by conforming finite elements in space and by the backward differentiation formula of order $k$ in time, with variable step-sizes for $k \le 3$ and constant step-size for $k = 4$.
The key ingredient of the analysis is a polynomial $(k+1)$-point time reconstruction of degree $k$ of the approximate solution combined with the elliptic reconstruction.
We conclude upper \emph{a~posteriori} error bounds in $L^\infty(L^2) \cap L^2(H^1)$ and in $L^\infty(H^1) \cap H^1(L^2)$.
Numerical experiments for $k = 1, \dots, 4$ with linear and quadratic finite elements confirm the optimal orders of the estimators in time and in space.
\end{abstract}


\maketitle
\section{Introduction}\label{SEC:introduction}

In this work, we derive reliable \emph{a~posteriori} error bounds for time discretization by the $k$-step backward differentiation formulas (\BDF{k}) for finite element discretizations of parabolic problems, using the heat equation as a model problem. 
In particular, we consider variable step-size \BDF{k} (\VSBDF{k}) for $k \le 3$ and uniform \BDF{4}. 
For an introduction on the \emph{a~posteriori} error analysis for parabolic partial differential equations, we refer to \cite{Smears:2025}.
Early \emph{a~posteriori} error estimates and adaptive algorithms for parabolic problems were derived in \cite{ErikssonJohnson:1991, Picasso:1998, Verfuerth:2003}, and upper bounds based on potential and flux reconstruction in \cite{ErnVohralik:2010, ErnSmearsVohralik:2017}.
For residual-based \emph{a~posteriori} error estimates of space-time finite element discretizations of parabolic equations, we refer to~\cite[Chap.~6]{Verfuerth:2013}.
For finite element discretizations, a common tool to obtain optimal-order \emph{a~posteriori} bounds is the elliptic reconstruction~\cite{MakridakisNochetto:2003,LakkisMakridakis:2006}, which reduces the spatial part of the error to \emph{a~posteriori} estimates for elliptic problems.
The key ingredient for the time discretization is a time reconstruction, whose residual is of optimal order.
Such reconstructions were developed for the Crank--Nicolson method~\cite{AkrivisMakridakisNochetto:2006, BaenschKarakatsaniMakridakis:2012, LozinskiET:2009}, for discontinuous Galerkin and implicit Runge--Kutta--Radau methods~\cite{MakridakisNochetto:2006}, and for Runge--Kutta collocation and continuous Galerkin methods of arbitrary order~\cite{AkrivisMakridakisNochetto:2009,AkrivisMakridakisNochetto:2011}.
Lower bounds within the elliptic reconstruction framework were established in~\cite{GeorgoulisMakridakis:2023}. 

Variable step-size BDF methods are widely used in adaptive solvers for stiff problems, see, e.g.,~\cite{HairerWanner:1996, ShampineReichelt:1997}.
Their stability on variable time grids requires restrictions on the ratio of consecutive time-steps \cite{Grigorieff:1983, CrouzeixLisbona:1984, CalvoGrandeGrigorieff:1990}, which were relaxed for \VSBDF{2} and \VSBDF{3} in \cite{LiLiao:2022}. 
For parabolic problems, the stability and convergence of \VSBDF{2} were analyzed,~e.g.,~in \cite{Becker:1998, Emmrich:2005, LiaoZhang:2021}, and of \VSBDF{3} in \cite{LiaoTangZhou:2023}.
For uniform time-steps, the energy technique based on the multipliers of \cite{NevanlinnaOdeh:1981} yields stability estimates for \BDF{k} with $k \le 5$~\cite{LubichET:2013, AkrivisLubich:2015}, and it was extended to \BDF{6} in \cite{AkrivisET:2021}.

For higher-order BDF methods, \emph{a~posteriori} error estimates are less developed.
Akrivis and Chatzipantelidis~\cite{AkrivisChatzipantelidis:2010} derived optimal-order estimates for the \BDF{2} method with constant time-steps.
Among their reconstructions is the three-point reconstruction,~i.e.,~the quadratic interpolant of three consecutive approximations, which was proposed in~\cite{LozinskiET:2009} for the Crank--Nicolson method.
Optimal estimators for the \BDF{2} method with variable step-sizes were derived in~\cite{WangMaoHuang:2022}, and \BDF{2}-based estimates were, for instance, obtained for the time-dependent Stokes equations~\cite{BaenschBrenner:2019}, and the Landau--Lifshitz--Gilbert equation~\cite{Karch:2026}.
To the best of our knowledge, reliable \emph{a~posteriori} error bounds for BDF methods of order $k \ge 3$ applied to parabolic problems are not available.

Our main contributions are the following.
We extend the three-point reconstruction to the general $(k+1)$-point time reconstruction $\recM{k}$, which interpolates $u_h^n, \dots, u_h^{n-k}$ by a polynomial of degree $k$ on each interval $(t_{n-1}, t_n]$.
We show that the residual of $\recM{k}$ consists of the interpolation error of the data and of terms involving the divided differences $\nabla^i u_h^n$ of order $i \ge k+1$.
Combined with the elliptic reconstruction, we conclude computable upper bounds for the error in $L^\infty(t_{k-1}, T; L^2(\Omega)) \cap L^2(t_{k-1}, T; H_0^1(\Omega))$ and in $L^\infty(t_{k-1}, T; H_0^1(\Omega)) \cap H^1(t_{k-1}, T; L^2(\Omega))$.
For $k \le 3$, the bounds hold for arbitrary sequences of time-steps.
Numerical experiments for $k = 1, \dots, 4$ with linear and quadratic finite elements confirm that the estimators converge with the optimal orders $k$ in time and $p$ in space, with effectivity indices that are essentially independent of $\tau$ and $h$.

The paper is organized as follows.
In Section~\ref{SEC:preliminaries}, we introduce the fully discrete scheme and the $(k+1)$-point time reconstruction.
Section~\ref{SEC:apost} contains the parabolic error equation, the definition of the estimators, and the \emph{a~posteriori} error bounds.
The residual identities and the bounds for the individual estimators are shown in Section~\ref{SEC:remaining-terms}.
Finally, Section~\ref{SEC:numerics} presents the numerical experiments.

\def\Sest{\SS}

\section{Preliminaries}\label{SEC:preliminaries}

In this section, we introduce the notation, weak formulation of the heat equation and its full discretization by finite elements in space and $\BDF{k}$ in time, as well as the $(k+1)$-point time reconstruction.

\subsection{Notation}
Let $\Omega \subset \R^d$, $d \in \{1, 2, 3\}$, be a bounded Lipschitz domain with polygonal ($d=2$) or polyhedral ($d=3$) boundary, and denote by $T > 0$ the final time.
We use the standard notation for Lebesgue and Sobolev spaces and denote by $H^{-1}(\Omega)$ the dual of $V \coloneqq H_0^1(\Omega)$. 
We denote by $\norm{\cdot}$ the norm of $L^2(\Omega)$, $\bL^2(\Omega) \coloneqq L^2(\Omega)^d$ and by $\dual{\cdot, \cdot}$ the $L^2(\Omega)$ inner product as well as the duality pairing between $H^{-1}(\Omega)$ and $H_0^1(\Omega)$.
For $v \in H_0^1(\Omega)$, we take $\norm{v}_{H^1(\Omega)} \coloneqq \norm{\grad v}_{L^2(\Omega)}$, and for $\star \in \{H^{-1}, L^2, H^1\}$ we abbreviate $\norm{\cdot}_\star \coloneqq \norm{\cdot}_{\star(\Omega)}$.
Further, let $a(v, w) \coloneqq \dual{\grad v, \grad w}$ for $v, w \in H_0^1(\Omega)$.
For a Banach space $X$ and an interval $I \subset \R$, we use the usual Bochner spaces $L^2(I; X)$, $L^\infty(I; X)$, $H^1(I; X)$ and $C(I; X)$.
Note that the bold symbol $\grad$ denotes the spatial gradient, whereas $\nabla^j$ with a superscript $j \ge 0$ denotes the backward divided differences in time defined later in \eqref{DEF:backward-div-diff}.

Let $\TT_h$ be a conforming and shape-regular triangulation of $\Omega$ into simplices $K$. 
Let $V_h \subset V$ be the finite element space of continuous piecewise polynomials of degree $p \ge 1$ on $\TT_h$.
Denote by $\Pi_h: L^2(\Omega) \to V_h$ the $L^2$-projection onto $V_h$.
For simplicity, we assume that the triangulation $\TT_h$ remains constant throughout the simulation.

In time, we consider a possibly non-uniform partition $0 = t_0 < t_1 < \dots < t_N = T$ with local step-sizes $\tau_n = t_n - t_{n-1}$ for $n = 1, \dots, N$, and set $I_n \coloneqq (t_{n-1}, t_n]$.
For a uniform partition, we write $\tau \coloneqq \tau_n = T / N$ for some $N \in \N$.
Moreover, we use the short-hand notation
\begin{align}
    \label{EQ:def-cm}
    c_m^n \coloneqq t_n - t_{n-m} = \sum_{j=0}^{m-1} \tau_{n-j}
\end{align}
for the sum of $m$ consecutive time-steps. 
Let $\P_j$ denote the polynomials of degree at most $j$, and $\P_j(X)$ the polynomials in $t$ of degree at most $j$ with coefficients in a vector space $X$.

For a sequence $\{ u_h^n \}_{n=0}^N \subset V_h$, we define the backward divided differences recursively for $n \ge j \ge 1$ by
\begin{align}
    \label{DEF:backward-div-diff}
    u_h[t_n,\ldots,t_{n-j}] 
    \coloneqq
    \nabla^j u_h^n 
    \coloneqq
    \frac{\nabla^{j-1}u_h^n-\nabla^{j-1}u_h^{n-1}}{t_n-t_{n-j}}, 
    \qq
    \nabla^0 u_h^n \coloneqq u_h^n.
\end{align}
We use the same notation for other sequences.
The corresponding Newton interpolation polynomial $P$ for the history $u_h^{n}, \ldots, u_h^{n-j}$ is given by
\begin{align}
    \label{EQ:Newton-interpol-poly}
    P[u_h^{n}, \ldots, u_h^{n-j}](t)
    =
    \sum_{i=0}^{j}
    \omega_i^{n}(t) 
    \nabla^i u_h^{n}
    \in \P_j(V_h),
    \qquad \omega_i^n(t) = \prod_{m=0}^{i-1}(t-t_{n-m}),
\end{align}
where an empty product evaluates to $1$.
The polynomial $P[u_h^{n}, \ldots, u_h^{n-j}](t)$ is the unique polynomial in $\P_j(V_h)$ that interpolates $u_h^{n-m}$ at $t_{n-m}$ for $m = 0, \dots, j$.

\subsection{Weak formulation and full discretization}
As model problem, we consider the heat equation with homogeneous Dirichlet boundary conditions
\begin{equation*}
\begin{aligned}
    \ptl_t u - \Delta u &= f \qq &&\textrm{in } \Omega \times (0, T), \\
    u &= 0 \qq &&\textrm{on } \ptl \Omega \times (0, T)
\end{aligned}
\end{equation*}
with initial data $u(0) = u_0 \in V$ and right-hand side $f \in C([0, T]; L^2(\Omega))$.
Its weak formulation is given by: Find $u \in L^2(0, T; V)$ with $\ptl_t u \in L^2(0, T; H^{-1}(\Omega))$ and $u(0) = u_0$ such that
\begin{align}
    \label{EQ:heateq}
    \dual{\ptl_t u, \varphi} + a(u, \varphi) = \dual{f, \varphi}, \qq \forall \varphi \in V
\end{align}
a.e.\ in $(0, T)$.
It is well known that the problem \eqref{EQ:heateq} admits a unique solution,~cf.~\cite[Sect.~7.1]{Evans:1998}.
Let $u_h^n$ denote the numerical approximation at time $t_n$.
Using the interpolation polynomial \eqref{EQ:Newton-interpol-poly}, the variable-step-size \BDF{k} (\VSBDF{k}) approximation of the time derivative at $t_n$, is given by 
\begin{align}
    \label{EQ:Def-vhnk}
    v_h^{n, (k)} \coloneqq \ptl_t P[u_h^{n}, \ldots, u_h^{n-k}](t_n)
    =
    \sum_{i=1}^{k}
    \ptl_t \omega_i^n(t_n) 
    \nabla^i u_h^{n}
\end{align}
for $n \ge k$,~cf.~\cite[Chap.~V]{HairerWanner:1996}.
Thus, we obtain the fully discrete scheme: 
Given starting values $u_h^0, \dots, u_h^{k-1} \in V_h$, find $u_h^n \in V_h$, $n = k, \dots, N$, such that
\begin{align}
    \label{EQ:heateq-discretized}
    \dual{v_h^{n, (k)}, \varphi_h} + a(u_h^n, \varphi_h) = \dual{f_h^n, \varphi_h}, \qq \forall \varphi_h \in V_h,
\end{align}
where $f_h^n \coloneqq \Pi_h f(t_n)$ for $n = 0, \dots, N$.

\subsection{Higher-order time reconstruction}

The key ingredient of our analysis is a reconstruction of the finite element approximations in time, whose residual is of optimal order.
For $n \ge k$, we define the local $(k+1)$-point time reconstruction as
\begin{align}\label{EQ:k+1-time-rec-def}
    \recMn{n}{k}(t) \coloneqq P[u_h^{n}, \ldots, u_h^{n-k}](t) \in \P_k(V_h),
\end{align}
which interpolates $u_h^{n-k}, \dots, u_h^n$ at $t_{n-k}, \dots, t_n$.
With $J \coloneqq [t_{k-1}, t_N]$, the (global) $(k+1)$-point time reconstruction $\recM{k} : J \to V_h$ is defined piecewise by
\begin{align}\label{EQ:k+1-time-rec-def-2}
    \recM{k}(t) \coloneqq \recMn{n}{k}(t) \qq \for t \in I_n,~n \ge k, 
\end{align}
and $\recM{k}(t_{k-1}) \coloneqq u_h^{k-1}$.
Since $\recMn{n}{k}(t_{n-1}) = u_h^{n-1}$ and $\recMn{n}{k}(t_n) = u_h^n$, the function $\recM{k}$ is continuous and piecewise polynomial in time. 
In particular, $\recM{k} \in H^1(J; V_h)$.
For $k=1$, we have the common piecewise linear (in time) reconstruction.
For $k=2$, we obtain the three-point reconstruction \cite{LozinskiET:2009, AkrivisChatzipantelidis:2010}.
Analogously to \eqref{EQ:k+1-time-rec-def}, we define the $(k+1)$-point reconstruction for the right-hand side by
\begin{align}\label{EQ:rhs-rec}
    \FF^{(k)}(t) = \sum_{i=0}^{k} \omega_i^{n}(t) \nabla^i f_h^n, \qq \for t \in I_n,~n \ge k.
\end{align}
The following lemma shows that the time derivative of $\recMn{n}{k}$ at $t_n$ coincides with the BDF approximation $v_h^{n, (k)}$ of the time derivative and provides the deviation within $I_n$.

\begin{Lemma}
    Let $\recMn{n}{k}$ be the $(k+1)$-point time reconstruction defined in \eqref{EQ:k+1-time-rec-def}. Then, it holds
    \begin{align*}
        \ptl_t \recMn{n}{k}(t_n) = v_h^{n, (k)}
    \end{align*}
    and
    \begin{align}
        \label{EQ:ptlt-recM}
        \partial_t\recMn{n}{k}(t)
        =
        v_h^{n,(k)}
        +
        \sum_{i=2}^{k}\rho_i^n(t)\,\nabla^iu_h^n 
    \end{align}
    for $t \in I_n$, $n \ge k$, where $\rho_i^n \in \mathbb{P}_{i-1}$ is defined for $i \ge 2$ as
    \begin{align}
        \label{EQ:Def_rho}
        \rho_i^n(t) \coloneqq \ptl_t \omega_i^n(t) - \ptl_t \omega_i^n(t_n),
    \end{align}
    where $\omega_i^{n}$ is defined in \eqref{EQ:Newton-interpol-poly}.
    In particular, $\ptl_t \recM{k}(t) = \ptl_t \recMn{n}{k}(t)$ for $t \in (t_{n-1}, t_n)$.
\end{Lemma}
\begin{proof}
    By the definition of $v_h^{n, (k)}$ in \eqref{EQ:Def-vhnk}, we have $\ptl_t \recMn{n}{k}(t_n) = v_h^{n, (k)}$.
    By differentiating the reconstruction, we obtain 
    \begin{align*}
        \partial_t\recMn{n}{k}(t)
        =
        \sum_{i=0}^{k}
        \ptl_t \omega_i^n(t) 
        \nabla^i u_h^{n}
        = 
        v_h^{n, (k)}
        +
        \sum_{i=0}^{k}
        (\ptl_t \omega_i^n(t) - \ptl_t \omega_i^n(t_n))
        \nabla^i u_h^{n},
    \end{align*}
    and with $\ptl_t \omega_0^n = 0$ and $\ptl_t \omega_1^n = 1$, we conclude \eqref{EQ:ptlt-recM}.
\end{proof}

Later, we require explicit formulas for the remainder polynomials $\rho_i^n$ introduced in \eqref{EQ:Def_rho}, which we collect in the following lemma.
The proof is provided in Appendix~\ref{AP:remaider-poly}.

\begin{Lemma}\label{LEM:rho_i}
    Let $\rho_i^n$, $i \ge 2$, be defined by \eqref{EQ:Def_rho}.
    Then, with $c_m^n$ defined in \eqref{EQ:def-cm}, we have
    \begin{equation}
    \label{EQ:rho-closed-form}
    \begin{aligned}
        \rho_2^n(t) &= 2(t-t_n),
        \\
        \rho_3^n(t) &= 3(t-t_n)^2 + 2(c_1^n + c_2^n)(t-t_n),
        \\ 
        \rho_4^n(t) &= 4(t-t_n)^3 + 3(c_1^n + c_2^n + c_3^n)(t-t_n)^2 + 2(c_1^n c_2^n + c_1^n c_3^n + c_2^n c_3^n)(t-t_n).
    \end{aligned}
    \end{equation}
\end{Lemma}

\section{\emph{A~posteriori} error estimates}\label{SEC:apost}

In this section, we first recall the elliptic reconstruction and residual-based \emph{a~posteriori} error estimates for elliptic problems. 
Then, we derive a parabolic error equation before introducing the \emph{a~posteriori} error estimators.
Finally, we provide the \emph{a~posteriori} error bounds.

\subsection{Elliptic reconstruction}
To split the error into temporal and spatial part, we use the elliptic reconstruction, which is based on the following discrete Laplacian.

\begin{Definition}
    \label{DEF: Discrete laplacoan}
    Let $v \in H_0^1(\Omega)$.
    Then the discrete Laplacian $- \Delta_h v \in V_h \subset H_0^1(\Omega)$ is defined by 
    \begin{align}
        \label{EQ: discrete Laplacian def}
        \dual{- \Delta_h v, \varphi_h} = a(v, \varphi_h), \qq \forall \varphi_h \in V_h.
    \end{align}
\end{Definition}

Using the discrete Laplacian, we introduce the elliptic reconstruction from \cite{LakkisMakridakis:2006, MakridakisNochetto:2003} as follows.

\begin{Definition}
    \label{DEF: Elliptic reconstruction}
    We define for $u_h \in V_h \subset H_0^1(\Omega)$ the elliptic reconstruction \hbox{$\RR u_h \in H_0^1(\Omega)$} by
    \begin{align}
        \label{EQ:ell-reconstruction}
        a({\RR u_h, \varphi}) = \dual{ - \Delta_h u_h, \varphi}, \qq \forall \varphi \in H_0^1(\Omega).
    \end{align}
\end{Definition}

Since $a(\RR v_h, \varphi_h) = \dual{- \Delta_h v_h, \varphi_h} = a(v_h, \varphi_h)$ for all $v_h, \varphi_h \in V_h$, the function $v_h$ is the Galerkin approximation of the solution $\RR v_h \in H_0^1(\Omega)$ of an elliptic problem with right-hand side $-\Delta_h v_h$.
Thus, $(\RR - I) v_h$ can be bounded by any \emph{a~posteriori} error estimate for elliptic problems.

\def\Cdual{\bar C_{-1}}
\def\CLtwo{\bar C_{0}}
\def\CHone{\bar C_{1}}

For $w_h \in V_h$, we define the element residual
\begin{align*}
    R_K(w_h) \coloneqq (\Delta - \Delta_h)w_h|_K
\end{align*}
for $K \in \TT_h$, and the face residual
\begin{align*}
    R_E(w_h) \coloneqq \jump{\ptl_{\bn} w_h}
\end{align*}
for $E \in \Sigma_h$, where $\Sigma_h$ denotes the set of interior faces of $\TT_h$ and $\jump{\ptl_{\bn} w_h}$ the jump of the normal derivative of $w_h$ across $E$ with respect to the unit normal vector $\bn$ of $E$.
For faces on $\ptl \Omega$, we set $R_E(w_h) = 0$.
For $K \in \TT_h$ and $E \in \Sigma_h$, let $h_K$ and $h_E$ denote the diameters of $K$ and $E$, respectively. 
These residuals define the standard residual-based \emph{a~posteriori} error estimators for elliptic problems, see,~e.g.,~\cite{Verfuerth:2013}.

We assume that the elliptic reconstruction satisfies the following residual-based \emph{a~posteriori} error bounds. 
Note that any other reliable elliptic \emph{a~posteriori} error estimator can be applied instead. 

\begin{Assumption}\label{ASS:a-post-bound-ell-rec}
    Let $p \ge 1$ be the polynomial degree of the finite element space $V_h$ and define $q \coloneqq \min\{p, 2\}$.
    We assume that for any $v_h \in V_h$, it holds
    \begin{align}
        \label{EQ:ell-rec-apost-dual}
        \norm{(\RR - I) v_h}_{H^{-1}}^2 &\le \eta(v_h; H^{-1}),
        \\
        \label{EQ:ell-rec-apost-L2}
        \norm{(\RR - I) v_h}^2 &\le \eta(v_h; L^2),
        \\
        \label{EQ:ell-rec-apost-H1}
        \norm{\grad (\RR - I) v_h}^2 &\le \eta(v_h; H^1),
    \end{align}
    where $\eta$ is defined by
    \begin{equation}
    \label{EQ: eta estimator standard}
   	\begin{aligned}
        \eta(v_h; H^{-1}) &\coloneqq \Cdual \Big( \sum\limits_{K \in \TT_h} h_K^{2q+2} \norm{R_K(v_h)}_{L^2(K)}^2 + \sum\limits_{E \in \Sigma_h} h_E^{2q+1} \norm{R_E(v_h)}_{L^2(E)}^2\Big),
        \\
   		\eta(v_h;  L^2) &\coloneqq \CLtwo \Big( \sum\limits_{K \in \TT_h} h_K^4 \norm{R_K(v_h)}_{L^2(K)}^2 + \sum\limits_{E \in \Sigma_h} h_E^3 \norm{R_E(v_h)}_{L^2(E)}^2 \Big),
   		\\
   		\eta(v_h; H^1) &\coloneqq \CHone \Big( \sum\limits_{K \in \TT_h} h_K^2 \norm{R_K(v_h)}_{L^2(K)}^2 + \sum\limits_{E \in \Sigma_h} h_E \norm{R_E(v_h)}_{L^2(E)}^2 \Big)
   	\end{aligned}
    \end{equation}
    for some constants $\Cdual, \CLtwo, \CHone > 0$.
\end{Assumption}

The bound \eqref{EQ:ell-rec-apost-H1} is the standard reliability estimate of the residual estimator,~cf.~\cite{Verfuerth:2013}, while \eqref{EQ:ell-rec-apost-dual} and \eqref{EQ:ell-rec-apost-L2} follow by duality arguments, provided that the dual problems are sufficiently regular.
With these preliminaries in place, we now derive the parabolic error equation on which both \emph{a~posteriori} error bounds are based.

\subsection{Parabolic error equation}

Let $\recW{k}$ be the interpolation polynomial of the pointwise (in time) elliptic reconstructions $\{\RR u_h^n\}_n$ according to \eqref{EQ:k+1-time-rec-def},~i.e.,~$\recWn{n}{k} \coloneqq P[\RR u_h^n, \dots, \RR u_h^{n-k}]$ for $n \ge k$, and $\recW{k}|_{I_n} \coloneqq \recWn{n}{k}$ as in \eqref{EQ:k+1-time-rec-def-2}.
Since $\RR$ is linear, we have $\recW{k} = \RR \recM{k}$.
Then, for all $\varphi \in V$, we have
\begin{align*}
    \dual{\ptl_t (\recM{k} - u), \varphi} + a(\recW{k-1} - u, \varphi)
    = 
    \dual{r_h^{(k)}, \varphi},
\end{align*}
where, with \eqref{EQ:ell-reconstruction}, the residual is given by 
\begin{equation}
\label{EQ:def-residual}
\begin{aligned}
    \dual{r_h^{(k)}, \varphi} &= - \dual{f, \varphi} 
    + \dual{\ptl_t \recM{k}, \varphi} + a(\recW{k-1}, \varphi)
    \\ &=
    \dual{\ptl_t \recM{k} - \Delta_h \recM{k-1} - f, \varphi}.
\end{aligned}
\end{equation}
With the error variables
\begin{align}\label{EQ:err-vars}
    \errM{k} = \recM{k} - u,\qq \errW{k} = \recW{k} - u,
\end{align}
we obtain the following parabolic error equation.

\begin{Lemma}
    For the solution $u$ of \eqref{EQ:heateq} and $u_h^n$ of     \eqref{EQ:heateq-discretized}, it holds for $t \in I_n$, $n \ge k$,
    \begin{align}
        \label{EQ:err-eq}
        \dual{\ptl_t \errM{k}, \varphi} + a(\errW{k-1}, \varphi) = \dual{r_h^{(k)}, \varphi}
    \end{align}
    for all $\varphi \in V$, where the residual $r_h^{(k)}$ is defined in \eqref{EQ:def-residual}.
\end{Lemma}
Introducing the elliptic-reconstruction defect 
\begin{align}\label{EQ:ell-rec-defect}
    \delta^{(k)} \coloneqq \recW{k} - \recM{k} = (\RR - I) \recM{k},
\end{align}
the error equation can be equivalently rewritten in the following formulation.

\begin{Corollary}
With $\delta^{(k)} = \recW{k} - \recM{k}$, the error equation \eqref{EQ:err-eq} is equivalent to
\begin{align}\label{EQ:err-eq-W}
    \dual{\ptl_t \errW{k}, \varphi} + a(\errW{k-1}, \varphi)
    = \dual{\ptl_t \delta^{(k)}, \varphi} + \dual{r_h^{(k)}, \varphi}, \qquad \forall \varphi \in V.
\end{align}
\end{Corollary}
\begin{proof}
    Follows immediately from $\errM{k} = \errW{k} - \delta^{(k)}$ and \eqref{EQ:err-eq}.
\end{proof}
The error equation \eqref{EQ:err-eq-W} is the starting point for the \emph{a~posteriori} theorems below.

\subsection{\emph{A~posteriori} error bounds}

To conveniently rewrite the residual \eqref{EQ:def-residual} later in Section~\ref{SCT:est-residual} for $k \le n \le 2k-2$, we introduce values $u_h^{-\ell}$ for $\ell = 1, \dots, k-1$ at negative times $t_{-\ell} = t_0 - \tau_1 \ell$.
Using these, the divided differences $\nabla^i u_h^n$ and the coefficients $c_m^n$ in \eqref{EQ:def-cm} are also defined for $n < i$ and $n < m$, respectively.
First of all, we define
\begin{align}
        \label{EQ:vhn-init-def}
        v_h^{n,(k)} \coloneqq \Delta_h u_h^n + f_h^n
\end{align}
for $1 \le n \le k-1$.
Then, we define $u_h^{-\ell}$ recursively for $\ell=1$ to $k-1$ such that 
\begin{align}\label{EQ:def-neg-time-vals-u}
    v_h^{n, (k)} = \Delta_h u_h^n + f_h^n = \sum\limits_{i=1}^k \ptl_t \omega_i^n(t_n) \nabla^i u_h^n
\end{align}
holds for $1 \le n \le k-1$. 
More details are given in Appendix~\ref{APP:approx-at-ge-times}.

This allows us to introduce the \emph{a~posteriori} error estimators derived in our analysis.

\begin{Definition}\label{DEF:collecting-estimators}
Let $\star \in \{H^{-1}, L^2, H^1\}$ denote the norm-type.
We define the \textbf{data approximation error estimator} by
\begin{align}\label{EQ:def-data-approx-err-est}
    \FF^{(k)}(\star) \coloneqq \int_J\norm{\FF^{(k-1)}-f}_\star^2\dt,
\end{align}
the \textbf{temporal truncation error estimator} for variable time-step \BDF{k}, $k\in\{1,2,3\}$, by 
\begin{equation}\label{EQ:def-NN1-NN2-NN3}
\begin{aligned}
    \NN^{(1)}(\star)&\coloneqq0,
    \qq
    \NN^{(2)}(\star) \coloneqq \sum_{n=2}^{N} \frac{(c_1^{n-1})^2(c_3^n)^2}{3}\,c_1^n\,
    \norm{\nabla^3u_h^n}_\star^2,
    \\
    \NN^{(3)}(\star) &\coloneqq \sum_{n=3}^{N}\Big[\Big(\frac{\tau_n^3}{3}\alpha_n^2
    + \frac{\tau_n^4}{6}\alpha_n\beta_n + \frac{\tau_n^5}{30}\beta_n^2\Big)
    \norm{\nabla^4u_h^n}_\star^2
    + \frac{\tau_n^5}{30}\gamma_n^2\,\norm{\nabla^5u_h^n}_\star^2\Big],
\end{aligned}
\end{equation}
where
\begin{align*}
    \alpha_n \coloneqq -\frac{c_1^{n-1}c_2^{n-1}c_4^n}{c_1^n},
    \qq
    \beta_n \coloneqq \frac{c_1^{n-2}c_2^{n-2}c_4^{n-1}}{c_1^{n-1}c_2^n}
        - \Big(\frac{c_1^{n-2}}{c_2^n}+\frac{c_2^{n-1}}{c_1^n}\Big)c_4^n,
    \qq
    \gamma_n \coloneqq -\frac{c_1^{n-2}c_2^{n-2}c_4^{n-1}c_5^n}{c_1^{n-1}c_2^n},
\end{align*}
and $c_j^n$ is defined in \eqref{EQ:def-cm}. For uniform \BDF{4}, we define
\begin{align}\label{EQ:def-NN4-uni}
    \NN^{(4)}(\star) &\coloneqq \sum_{n=4}^N\Big[\frac{10730}{21}\tau^9\norm{\nabla^5u_h^n}_\star^2
    + \frac{15790}{21}\tau^{11}\norm{\nabla^6u_h^n}_\star^2
    + 3360\,\tau^{13}\norm{\nabla^7u_h^n}_\star^2\Big].
\end{align}
Moreover, we define the \textbf{elliptic reconstruction error estimator}
$\Theta^{(k)}(\star)$, and its $L^\infty$ in-time
form $\Theta_\infty^{(k)}(\star)$, as the fully computable bounds satisfying
\begin{align*}
    \int_J\norm{\delta^{(k)}}_\star^2\dt \le \Theta^{(k)}(\star),
    \qq
    \norm{\delta^{(k)}}_{L^\infty(J;\star)}^2 \le \Theta_\infty^{(k)}(\star).
\end{align*}
The explicit definitions are given in Lemma~\ref{LEM:grad-delta-est} and
Lemma~\ref{LEM:delta-Linfty-est}. 
Similarly, we define the \textbf{reconstruction-rate error estimator} $\Lambda^{(k)}(\star)$ by
\begin{align*}
    \int_J\norm{\ptl_t\delta^{(k)}}_\star^2\dt \le \Lambda^{(k)}(\star)
\end{align*}
with explicit definition given in Lemma~\ref{LEM:dt-delta-est}.
Further, we define the \textbf{inter-order error estimator} by
\begin{equation}
\begin{aligned}
    \Sest^{(k)}(H^1) &\coloneqq \sum\limits_{n = k}^N \omega^{n, (k)} \Big(\norm{\nabla^k (\grad u_h^n)}^2 + \eta(\nabla^k u_h^n; H^1) \Big),
    \\
    \Sest^{(k)}(H^2) &\coloneqq \sum\limits_{n = k}^N \omega^{n, (k)} \norm{\nabla^k (\Delta_hu_h^n)}^2,
\end{aligned}
\end{equation}
where $\omega^{n, (k)} \coloneqq \int_{t_{n-1}}^{t_n} |\omega_k^n(t)|^2 \dt$, $n \ge k$,
and finally the \textbf{initial error estimator} by
\begin{align}\label{EQ:Def-init-err-est}
    \II(\star) \coloneqq \norm{u_h^{k-1}-u(t_{k-1})}_\star^2+\eta(u_h^{k-1};\star).
\end{align}
\end{Definition}

Additionally, we introduce 
\begin{align}
    \rho^{(k)} \coloneqq \recW{k-1} - \recW{k}, \qq \rho_h^{(k)} \coloneqq \recM{k-1} - \recM{k}
\end{align}
to split the error components. 
Note that $\rho^{(k)} = \RR \rho_h^{(k)}$ and $\rho_h^{(k)}(t) = - \omega_k^n(t) \nabla^k u_h^n$ for $t \in I_n$.
Next, we estimate the error \eqref{EQ:err-vars} in the $L^\infty(L^2)$- and $L^2(H^1)$-norm.

\begin{Theorem}[\emph{A~posteriori} error bound \rom{1}]\label{THM:first-apost-pointwise}
    Let $u$ be the solution of \eqref{EQ:heateq} and $u_h^n$ of \eqref{EQ:heateq-discretized} for $n \ge k$. Then, for every $s\in[t_{k-1},T]$ it holds for the error \eqref{EQ:err-vars}
    \begin{equation}
    \begin{aligned}
        \norm{\errM{k}(s)}^2 + \int_{t_{k-1}}^{s} \norm{\grad \errM{k}}^2 \dt
        &\le         
        2 \norm{\delta^{(k)}(s)}^2
        + 2\int_{t_{k-1}}^{s}\norm{\grad\delta^{(k)}}^2\dt
        \\&\qq
        + 4\norm{\errM{k}(t_{k-1})}^2
        + 4\norm{\delta^{(k)}(t_{k-1})}^2
        \\&\qq
        +6\int_{t_{k-1}}^{s} \Big( \norm{\ptl_t\delta^{(k)}}_{H^{-1}}^2+\norm{r_h^{(k)}}_{H^{-1}}^2+\norm{\grad\rho^{(k)}}^2 \Big)\dt
        .
    \end{aligned}
    \end{equation}
\end{Theorem}
\begin{proof}
    Testing \eqref{EQ:err-eq-W} with $\varphi = \errW{k}$ yields
    \begin{align*}
        \frac12 \ptl_t \norm{\errW{k}}^2 + a(\errW{k-1}, \errW{k}) = \dual{\ptl_t \delta^{(k)}, \errW{k}} + \dual{r_h^{(k)}, \errW{k}}.
    \end{align*}
    Using
    \begin{align*}
        a(\errW{k-1}, \errW{k}) &= \norm{\grad \errW{k}}^2 + a(\errW{k-1} - \errW{k}, \errW{k})
        = 
        \norm{\grad \errW{k}}^2 + \dual{\grad \rho^{(k)}, \grad \errW{k}}
    \end{align*}
    we obtain
    \begin{align*}
        \frac12 \ptl_t \norm{\errW{k}}^2 + \norm{\grad \errW{k}}^2 = \dual{\ptl_t \delta^{(k)}, \errW{k}} + \dual{r_h^{(k)}, \errW{k}} - \dual{\grad \rho^{(k)}, \grad \errW{k}}
    \end{align*}
    and applying Young's inequality gives
    \begin{align*}
        \ptl_t\norm{\errW{k}}^2+\norm{\grad\errW{k}}^2 \le 3\Big(\norm{\ptl_t\delta^{(k)}}_{H^{-1}}^2+\norm{r_h^{(k)}}_{H^{-1}}^2+\norm{\grad\rho^{(k)}}^2\Big) \eqqcolon 3F(t).
    \end{align*}
    Integrating from $t_{k-1}$ to $s$, for every $s\in[t_{k-1},T]$, yields
    \begin{align}\label{EQ:a-post-intermed}
        \norm{\errW{k}(s)}^2+\int_{t_{k-1}}^{s}\norm{\grad\errW{k}}^2 \dt \le \norm{\errW{k}(t_{k-1})}^2 + 3\int_{t_{k-1}}^{s} F \dt.
    \end{align}
    Using $\errW{k} = \errM{k} + \delta^{(k)}$, as well as $(a+b)^2 \le 2 a^2 + 2 b^2$, we obtain
    \begin{equation}\label{EQ:a-post-intermed2}
    \begin{aligned}
        \norm{\errM{k}(s)}^2+\int_{t_{k-1}}^{s}\norm{\grad\errM{k}}^2 \dt
        &\le 
        2\norm{\delta^{(k)}(s)}^2
        +2\int_{t_{k-1}}^{s}\norm{\grad\delta^{(k)}}^2\dt
        +6\int_{t_{k-1}}^{s}F\dt
        \\&\qq
        +4\norm{\errM{k}(t_{k-1})}^2
        +4\norm{\delta^{(k)}(t_{k-1})}^2
        ,
    \end{aligned}
    \end{equation}
    which concludes the proof.
\end{proof}

As an immediate consequence, taking the supremum and the $L^2(J)$-integral of \eqref{EQ:a-post-intermed}, and proceeding as for \eqref{EQ:a-post-intermed2}, we obtain the following combined bound.

\begin{Corollary}[$L^\infty$ \emph{a~posteriori} error bound \rom{1}]\label{COR:first-apost-Linfty}
    Let $u$ be the solution of \eqref{EQ:heateq} and $u_h^n$ of \eqref{EQ:heateq-discretized} for $n \ge k$.
    Then, it holds for the error \eqref{EQ:err-vars}
    \begin{equation}
    \begin{aligned}
        \norm{\errM{k}}_{L^\infty(J;L^2(\Omega))}^2 + \norm{\grad\errM{k}}_{L^2(J;\bL^2(\Omega))}^2
        \le
        8\II(L^2)
        &+ 12k\Big(\FF^{(k)}(H^{-1})+\NN^{(k)}(H^{-1})\Big)
        + 12\Lambda^{(k)}(H^{-1})
        \\&+ 12 \Sest^{(k)}(H^1)
        + 2\Theta_\infty^{(k)}(L^2)
        + 2\Theta^{(k)}(H^1).
    \end{aligned}
    \end{equation}
\end{Corollary}
\begin{proof}
    By \eqref{EQ:a-post-intermed}, we have
    \begin{align*}
        \norm{\errW{k}(s)}^2 \le \norm{\errW{k}(t_{k-1})}^2 + 3\int_{t_{k-1}}^{s} F \dt, 
        \qq
        \int_{t_{k-1}}^{s}\norm{\grad\errW{k}}^2 \dt \le \norm{\errW{k}(t_{k-1})}^2 + 3\int_{t_{k-1}}^{s} F \dt
    \end{align*}
    for $s \in J$, where $F = \norm{\ptl_t\delta^{(k)}}_{H^{-1}}^2+\norm{r_h^{(k)}}_{H^{-1}}^2+\norm{\grad\rho^{(k)}}^2$.
    Taking the supremum in the first and $s=T$ in the second inequality, then adding both, and proceeding as for \eqref{EQ:a-post-intermed2}, we obtain
    \begin{align*}
        \norm{\errM{k}}_{L^\infty(J; L^2(\Omega))}^2+\int_{t_{k-1}}^{T}\norm{\grad\errM{k}}^2 \dt
        &\le 
        2\norm{\delta^{(k)}}_{L^\infty(J; L^2(\Omega))}^2
        +2\int_{t_{k-1}}^{T}\norm{\grad\delta^{(k)}}^2\dt
        +12\int_{t_{k-1}}^{T}F\dt
        \\&\qq
        +8\norm{\errM{k}(t_{k-1})}^2
        +8\norm{\delta^{(k)}(t_{k-1})}^2
        .
    \end{align*}
    Finally, applying 
    Lemma~\ref{LEM:initial-error} ($\II$),
    Lemma~\ref{LEM:est-residual-final} ($\FF^{(k)},\NN^{(k)}$),
    Lemma~\ref{LEM:dt-delta-est} ($\Lambda^{(k)}$), 
    Lemma~\ref{LEM:grad-delta-est} ($\Theta^{(k)}$), 
    Lemma~\ref{LEM:delta-Linfty-est} ($\Theta_\infty^{(k)}$), and
    Lemma~\ref{LEM:est-rho} ($\Sest^{(k)}$) we conclude the proof.
\end{proof}

We now consider the second \emph{a~posteriori} bound, which controls $\ptl_t\errM{k}$ and $\grad\errM{k}$ instead, and is derived by testing \eqref{EQ:err-eq-W} with $\ptl_t\errW{k}$ instead of $\errW{k}$.

\begin{Theorem}[\emph{A~posteriori} error bound \rom{2}]\label{THM:second-apost-pointwise}
    Let $u$ be the solution of \eqref{EQ:heateq} and $u_h^n$ of \eqref{EQ:heateq-discretized} for $n \ge k$. Then, for every $s\in[t_{k-1},T]$ it holds for the error \eqref{EQ:err-vars}
    \begin{equation}
    \begin{aligned}
        \norm{\grad\errM{k}(s)}^2 + \int_{t_{k-1}}^{s} \norm{\ptl_t\errM{k}}^2 \dt
        &\le
        2\norm{\grad\delta^{(k)}(s)}^2
        + 8\int_{t_{k-1}}^{s}\norm{\ptl_t\delta^{(k)}}^2\dt
        + 2\norm{\grad\errW{k}(t_{k-1})}^2
        \\&\qq
        + 6\int_{t_{k-1}}^{s}\Big(\norm{r_h^{(k)}}^2+\norm{\Delta_h\rho_h^{(k)}}^2\Big)\dt
        .
    \end{aligned}
    \end{equation}
\end{Theorem}
\begin{proof}
    Testing \eqref{EQ:err-eq-W} with $\varphi = \ptl_t\errW{k}$ yields
    \begin{align*}
        \norm{\ptl_t \errW{k}}^2 + a(\errW{k-1}, \ptl_t\errW{k}) = \dual{\ptl_t \delta^{(k)}, \ptl_t\errW{k}} + \dual{r_h^{(k)}, \ptl_t\errW{k}}.
    \end{align*}
    Using \eqref{EQ:ell-reconstruction} and
    \begin{align*}
        a(\errW{k-1}, \ptl_t\errW{k}) &= \frac12 \ptl_t \norm{\grad \errW{k}}^2 + a(\errW{k-1} - \errW{k}, \ptl_t\errW{k})
        =
        \frac12 \ptl_t \norm{\grad \errW{k}}^2 - \dual{\Delta_h\rho_h^{(k)}, \ptl_t\errW{k}}
    \end{align*}
    we obtain
    \begin{align*}
        \norm{\ptl_t\errW{k}}^2 + \frac12\ptl_t \norm{\grad \errW{k}}^2 = \dual{\ptl_t \delta^{(k)}, \ptl_t\errW{k}} + \dual{r_h^{(k)}, \ptl_t\errW{k}} + \dual{\Delta_h\rho_h^{(k)}, \ptl_t\errW{k}}.
    \end{align*}
    Applying Young's inequality gives
    \begin{align*}
        \norm{\ptl_t\errW{k}}^2+\ptl_t\norm{\grad\errW{k}}^2 \le 3\Big(\norm{\ptl_t\delta^{(k)}}^2+\norm{r_h^{(k)}}^2+\norm{\Delta_h\rho_h^{(k)}}^2\Big) \eqqcolon 3F_2(t).
    \end{align*}
    Integrating from $t_{k-1}$ to $s$ yields, for every $s\in[t_{k-1},T]$,
    \begin{align}\label{EQ:second-apost-intermed}
        \int_{t_{k-1}}^{s}\norm{\ptl_t\errW{k}}^2\dt + \norm{\grad\errW{k}(s)}^2 \le \norm{\grad\errW{k}(t_{k-1})}^2 + 3\int_{t_{k-1}}^{s} F_2 \dt.
    \end{align}
    Using $\errW{k} = \errM{k}+\delta^{(k)}$, as well as $(a+b)^2 \le 2 a^2 + 2 b^2$, we obtain
    \begin{align*}
        \norm{\grad\errM{k}(s)}^2+\int_{t_{k-1}}^{s}\norm{\ptl_t\errM{k}}^2 \dt
        &\le
        2\norm{\grad\delta^{(k)}(s)}^2
        +2\int_{t_{k-1}}^{s}\norm{\ptl_t\delta^{(k)}}^2\dt
        \\&\qq
        +6\int_{t_{k-1}}^{s}F_2\dt
        +2\norm{\grad\errW{k}(t_{k-1})}^2
        .
    \end{align*}
    Since $F_2$ already contains $\norm{\ptl_t\delta^{(k)}}^2$, the two occurrences of
    $\int_{t_{k-1}}^{s}\norm{\ptl_t\delta^{(k)}}^2\dt$ combine, which concludes the proof.
\end{proof}

As before, the combined $L^\infty$/$L^2$-in-time formulation follows by taking the supremum and the $L^2(J)$-integral of \eqref{EQ:second-apost-intermed}, and repeating the argument used for Theorem~\ref{THM:first-apost-pointwise}.

\begin{Corollary}[$L^\infty$ \emph{a~posteriori} error bound \rom{2}]\label{COR:second-apost-Linfty}
    Let $u$ be the solution of \eqref{EQ:heateq} and $u_h^n$ of \eqref{EQ:heateq-discretized} for $n \ge k$.
    Then, it holds for the error \eqref{EQ:err-vars}
    \begin{equation}
    \begin{aligned}
        \norm{\grad\errM{k}}_{L^\infty(J;\bL^2(\Omega))}^2 + \norm{\ptl_t\errM{k}}_{L^2(J;L^2(\Omega))}^2
        \le
        8\II(H^1)
        &+ 12k\Big(\FF^{(k)}(L^2)+\NN^{(k)}(L^2)\Big)
        + 14\Lambda^{(k)}(L^2)
        \\&+ 12\Sest^{(k)}(H^2)
        + 2\Theta_\infty^{(k)}(H^1).
    \end{aligned}
    \end{equation}
\end{Corollary}
\begin{proof}
    By \eqref{EQ:second-apost-intermed}, we have
    \begin{align*}
        \norm{\grad\errW{k}(s)}^2 &\le \norm{\grad\errW{k}(t_{k-1})}^2 + 3\int_{t_{k-1}}^{s} F_2 \dt,
        \\
        \int_{t_{k-1}}^{s}\norm{\ptl_t\errW{k}}^2\dt &\le \norm{\grad\errW{k}(t_{k-1})}^2 + 3\int_{t_{k-1}}^{s} F_2 \dt
    \end{align*}
    for $s \in J$, where $F_2 = \norm{\ptl_t\delta^{(k)}}^2+\norm{r_h^{(k)}}^2+\norm{\Delta_h\rho_h^{(k)}}^2$.
    Taking the supremum in the first and $s=T$ in the second, then adding both, and proceeding as in the proof of Corollary~\ref{COR:first-apost-Linfty}, we obtain
    \begin{align*}
        \norm{\grad\errM{k}}_{L^\infty(J;\bL^2(\Omega))}^2+\int_{t_{k-1}}^{T}\norm{\ptl_t\errM{k}}^2 \dt
        &\le
        2\norm{\grad\delta^{(k)}}_{L^\infty(J; \bL^2(\Omega))}^2
        \\&\qq
        +14\int_{t_{k-1}}^{T}\norm{\ptl_t\delta^{(k)}}^2\dt
        +4\norm{\grad\errW{k}(t_{k-1})}^2
        \\&\qq
        +12\int_{t_{k-1}}^{T}\Big(\norm{r_h^{(k)}}^2+\norm{\Delta_h\rho_h^{(k)}}^2\Big)\dt
        .
    \end{align*}
    Finally, applying
    Lemma~\ref{LEM:initial-error} ($\II$),
    Lemma~\ref{LEM:est-residual-final} ($\FF^{(k)},\NN^{(k)}$),
    Lemma~\ref{LEM:dt-delta-est} ($\Lambda^{(k)}$),
    Lemma~\ref{LEM:est-rho} ($\Sest^{(k)}$), and
    Lemma~\ref{LEM:delta-Linfty-est} ($\Theta_\infty^{(k)}$) we conclude the proof.
\end{proof}

\section{Estimating the remaining terms}\label{SEC:remaining-terms}
It remains to establish the explicit bounds for each estimator introduced in Definition~\ref{DEF:collecting-estimators} and used in Corollaries~\ref{COR:first-apost-Linfty} and \ref{COR:second-apost-Linfty} above.
We begin with the initial error.
 
\begin{Lemma}\label{LEM:initial-error}
    Let $\star \in \{L^2, H^1\}$.
    For the initial error, it holds
    \begin{align*}
        \norm{\errM{k}(t_{k-1})}_\star^2 + \norm{\delta^{(k)}(t_{k-1})}_\star^2 \le \norm{u_h^{k-1} - u(t_{k-1})}_\star^2 + \eta(u_h^{k-1}; \star) = \II(\star).
    \end{align*}
\end{Lemma}
\begin{proof}
    By \eqref{EQ:k+1-time-rec-def} and \eqref{EQ:k+1-time-rec-def-2}, we have
    \begin{align*}
        \recM{k}(t_{k-1}) = \recMn{k}{k}(t_{k-1}) = u_h^k - c_1^k\nabla^1u_h^k
        = u_h^k - c_1^k\cdot\frac{u_h^k-u_h^{k-1}}{c_1^k} = u_h^{k-1}
    \end{align*}
    Since $\recW{k}$ is built by the identical formula with $u_h^n$ replaced by $\RR u_h^n$ throughout, we have $\recW{k}(t_{k-1}) = \RR u_h^{k-1}$. 
    Consequently,
    \begin{align*}
        \errM{k}(t_{k-1}) = \recM{k}(t_{k-1}) - u(t_{k-1}) = u_h^{k-1} - u(t_{k-1}),
        \qquad
        \delta^{(k)}(t_{k-1}) = \RR u_h^{k-1} - u_h^{k-1}.
    \end{align*}
    Applying \eqref{EQ:ell-rec-apost-L2} for $\star=L^2$ or \eqref{EQ:ell-rec-apost-H1} for $\star=H^1$ concludes the proof.
\end{proof}

\subsection{Estimating the residual}\label{SCT:est-residual}
In the following Lemma, we rewrite the residual in terms of higher-order backward differences and the right-hand side reconstruction for variable time-steps for \BDF{k}, $k \in \{1, 2, 3\}$.

\begin{Lemma}\label{LEM:residual-identities}
    Let $r_h^{(k)}$ be the residual given in \eqref{EQ:def-residual}.
    Then, it holds
    \begin{align*}
        r_h^{(1)}(t) &= \FF^{(0)}(t) - f(t),
        \\
        r_h^{(2)}(t) 
        &= 
        \FF^{(1)}(t) - f(t) - \frac{c_1^{n-1} c_3^n}{c_1^n} \omega_1^{n}(t) \nabla^3 u_h^n,
        \\
        r_h^{(3)}(t)
        &= \FF^{(2)}(t) - f(t)
        + A_4(t) \nabla^4 u_h^n
        + A_5(t) \nabla^5 u_h^n,
    \end{align*}
    for $t \in I_n$, for $n \ge k$, where the coefficients are given by
    \begin{equation}\label{EQ:Def-A4-A5}
    \begin{aligned}
        A_4(t) &= 
        - \frac{c_1^{n-1} c_2^{n-1} c_4^n }{c_1^n} \omega_1^{n}(t)
        + \left(\frac{c_1^{n-2} c_2^{n-2} c_4^{n-1}}{c_1^{n-1} c_2^n} - \bigg( \frac{c_1^{n-2}}{c_2^n} + \frac{c_2^{n-1}}{c_1^n} \bigg) c_4^{n}\right) \omega_2^{n}(t) ,
        \\
        A_5(t) &= - \frac{c_1^{n-2} c_2^{n-2} c_4^{n-1} c_5^n}{c_1^{n-1} c_2^n} \omega_2^{n}(t)
    \end{aligned}
    \end{equation}
    and the values $u_h^{-\ell}$ at negative times $t_{-\ell}$ for $\ell = 1, \dots, k-1$ are defined by \eqref{EQ:def-neg-time-vals-u}.
\end{Lemma}
\begin{proof}
    By \eqref{EQ:def-residual}, the residual is given by $r_h^{(k)}(t) = - f(t) + \ptl_t \recM{k}(t) - \Delta_h \recM{k-1}(t)$ for $t \in I_n$, $n \ge k$.
    Using \eqref{EQ:heateq-discretized} and \eqref{EQ: discrete Laplacian def}, we obtain 
    \begin{align}\label{EQ:pw-heat-eq}
        v_h^{n, (k)} - \Delta_h u_h^n - f_h^n = 0
    \end{align}
    for $n \ge k$ since 
    \begin{align*}
        \dual{v_h^{n, (k)} - \Delta_h u_h^n - f_h^n, \varphi} &=
        \dual{v_h^{n, (k)} - \Delta_h u_h^n - f_h^n, \Pi_h \varphi} 
        \\ &=
        \dual{v_h^{n, (k)}, \Pi_h \varphi} + a(u_h^n, \Pi_h\varphi) - \dual{f_h^n, \Pi_h \varphi} = 0.
    \end{align*}
    Thus, with \eqref{EQ:k+1-time-rec-def} and \eqref{EQ:ptlt-recM} we rewrite the residual by
    \begin{equation} \label{EQ-RR-sums}
    \begin{aligned}
        r_h^{(k)}(t) &= 
        v_h^{n,(k)} - f(t)
        + \sum_{i=2}^k \rho_i^n(t) \nabla^i u_h^n
        - \sum_{i=0}^{k-1} \omega_i^{n}(t) \Delta_h \big( \nabla^iu_h^n \big)
    \end{aligned}
    \end{equation}
    for $t \in I_n$, $n \ge k$.
    Next, define 
    \begin{align}\label{EQ:def-qhn}
        q_h^n \coloneqq \nabla^1 u_h^n - \Delta_h u_h^n - f_h^n
    \end{align}
    for $n \ge 1$.
    Using \eqref{EQ:pw-heat-eq} for $n \ge k$ and \eqref{EQ:vhn-init-def} for $1 \le n \le k-1$ we obtain 
    \begin{align}\label{EQ:qhn-eq}
        q_h^n = \nabla^1 u_h^n - \Delta_h u_h^n - f_h^n = \nabla^1 u_h^n - v_h^{n, (k)}
    \end{align}
    for $n \ge 1$, which yields $\Delta_h u_h^n = \nabla^1 u_h^n - q_h^n - f_h^n$ for $n \ge 1$.
    Substituting $\Delta_h (\nabla^i u_h^n) = \nabla^i (\Delta_h u_h^n) = \nabla^i(\nabla^1 u_h^n) - \nabla^i q_h^n - \nabla^i f_h^n$ for $n \ge k > i \ge 0$ into \eqref{EQ-RR-sums} and using \eqref{EQ:rhs-rec}, we obtain
    \begin{equation} \label{EQ-RR-master}
        r_h^{(k)}(t) = \FF^{(k-1)}(t) - f(t) + \sum_{i=1}^{k-1} \omega_i^{n}(t) \nabla^i q_h^n + \sum_{i=2}^k \rho_i^n(t) \nabla^i u_h^n - \sum_{i=1}^{k-1} \omega_i^{n}(t) \nabla^i \big( \nabla^1 u_h^n \big)
    \end{equation}
    for $n \ge k$.
    Further, using \eqref{EQ:Def-vhnk} and \eqref{EQ:def-neg-time-vals-u}, we have 
    \begin{align}\label{EQ:vhn-n-ge-1}
        v_h^{n, (k)} = \sum\limits_{i=1}^k \ptl_t \omega_i^n(t_n) \nabla^i u_h^n
    \end{align}
    for $n \ge 1$.
        
    \textbf{Case $k=1$:} The sums in \eqref{EQ-RR-master} vanish, leaving
    \begin{align*}
        r_h^{(1)}(t) &= 
        \FF^{(0)}(t) - f(t).
    \end{align*}

    \textbf{Case $k=2$:}
    We evaluate the divided difference of the first divided difference
    \begin{align*}
        \nabla^1 \big( \nabla^1 u_h^n \big) 
        = \frac{\nabla^1 u_h^n - \nabla^1 u_h^{n-1}}{c_1^n} 
        = \frac{c_2^n}{c_1^n} \nabla^2 u_h^n,
    \end{align*}
    where $c_i^n$ is defined in \eqref{EQ:def-cm}.
    Substituting this into \eqref{EQ-RR-master} for $k=2$ yields
    \begin{align}
        r_h^{(2)}(t)
        &= \FF^{(1)}(t) - f(t) + \omega_1^{n}(t) \nabla^1 q_h^n + \rho_2^n(t) \nabla^2 u_h^n - \omega_1^{n}(t) \frac{c_2^n}{c_1^n} \nabla^2 u_h^n
        \\&= \FF^{(1)}(t) - f(t) + \omega_1^{n}(t) \nabla^1 q_h^n + \omega_1^{n}(t) \left( 2 - \frac{c_2^n}{c_1^n} \right) \nabla^2 u_h^n. \label{EQ:R2-var1}
    \end{align}
    By \eqref{EQ:vhn-n-ge-1} and \eqref{EQ:def-qhn}, we have $q_h^n = \nabla^1 u_h^n - \ptl_t \omega_1^{n}(t_n) \nabla^1 u_h^n - \ptl_t \omega_2^{n}(t_n) \nabla^2 u_h^n 
    = -c_1^n\nabla^2u_h^n$.
    Applying the first-order divided difference yields
    \begin{align*}
        \nabla^1 q_h^n = \frac{q_h^n - q_h^{n-1}}{c_1^n} = -\nabla^2u_h^n + \frac{c_1^{n-1}}{c_1^n}\nabla^2u_h^{n-1}.
    \end{align*}
    Combining the terms on the right-hand side of \eqref{EQ:R2-var1}, we obtain 
    \begin{align*}
        \nabla^1q_h^n + \Big(2-\frac{c_2^n}{c_1^n}\Big)\nabla^2u_h^n = \frac{c_1^{n-1}}{c_1^n}\big[\nabla^2u_h^{n-1}-\nabla^2u_h^n\big] = -\frac{c_1^{n-1}}{c_1^n}
        c_3^n \nabla^3u_h^n.
    \end{align*}
    Thus, in total we conclude
    \begin{align*}
        r_h^{(2)}(t)
        &= \FF^{(1)}(t) - f(t) - \frac{c_1^{n-1}}{c_1^n}\, c_3^n\, \omega_1^{n}(t) \nabla^3 u_h^n.
    \end{align*}

    \textbf{Case $k=3$:}
    We require the second divided difference of $\nabla^1 u_h^n$. Using $\nabla^2 u_h^{n-1} = \nabla^2 u_h^n - c_3^n \nabla^3 u_h^n$, we find
    \begin{align*}
        \nabla^2 \big( \nabla^1 u_h^n \big)
        &= \frac{\nabla^1(\nabla^1 u_h^n) - \nabla^1(\nabla^1 u_h^{n-1})}{c_2^n}
        \\&= \frac{1}{c_2^n} \left( \frac{c_2^n}{c_1^n} \nabla^2 u_h^n - \frac{c_2^{n-1}}{c_1^{n-1}} \nabla^2 u_h^{n-1} \right)
        \\&= \frac{1}{c_2^n} \left( \frac{c_2^n}{c_1^n} - \frac{c_2^{n-1}}{c_1^{n-1}} \right) \nabla^2 u_h^n + \frac{c_2^{n-1} c_3^n}{c_1^{n-1} c_2^n} \nabla^3 u_h^n.
    \end{align*}
    Substituting this and the $k=2$ result into \eqref{EQ-RR-master} for $k=3$ groups the terms as
    \begin{align}
        \label{EQ:RR3-tmp1}
        r_h^{(3)}(t)
        &= \FF^{(2)}(t) - f(t) 
        + \omega_1^{n}(t) \nabla^1 q_h^n + \omega_2^{n}(t) \nabla^2 q_h^n
        + A_2(t) \nabla^2 u_h^n 
        + A_3(t) \nabla^3 u_h^n,
    \end{align}
    with the coefficients given by
    \begin{align*}
        A_2(t) &= \Big(2 - \frac{c_2^n}{c_1^n}\Big)\omega_1^{n}(t) - \frac{1}{c_2^n} \left( \frac{c_2^n}{c_1^n} - \frac{c_2^{n-1}}{c_1^{n-1}} \right) \omega_2^{n}(t),
        \\
        A_3(t) &= \rho_3^n(t) - \frac{c_2^{n-1} c_3^n}{c_1^{n-1} c_2^n} \omega_2^{n}(t) = 3 \omega_1^{n}(t)^2 + 2 (c_1^n + c_2^n) \omega_1^{n}(t) - \frac{c_2^{n-1} c_3^n}{c_1^{n-1} c_2^n} \omega_2^{n}(t)
        \\&=
        \omega_2^{n}(t) \Big( 3 - \frac{c_2^{n-1} c_3^n}{c_1^{n-1} c_2^n} \Big)
        + \big( 2 c_2^n - c_1^n \big) \omega_1^{n}(t).
    \end{align*}
    Further, we have with \eqref{EQ:vhn-n-ge-1}, \eqref{EQ:def-qhn},
    $\ptl_t \omega_2^{n}(t_n) = \ptl_t \omega_1^{n}(t_n) c_1^n = c_1^n$ 
    and $\ptl_t \omega_3^{n}(t_n) = \ptl_t \omega_2^{n}(t_n) c_2^n = c_2^n c_1^n$ that
    \begin{align}
        q_h^n &= 
        \nabla^1 u_h^n - \ptl_t \omega_1^{n}(t_n) \nabla^1 u_h^n 
        - \ptl_t \omega_2^{n}(t_n) \nabla^2 u_h^n
        - \ptl_t \omega_3^{n}(t_n) \nabla^3 u_h^n
        \\&=
        - c_1^n \nabla^2 u_h^n
        - c_1^n c_2^n \nabla^3 u_h^n.
    \end{align}
    Collecting the terms with $\omega_2^{n}$ and resolving $\nabla^2 q_h^n$ by 
    \begin{align*}
        \nabla^2 q_h^n = \frac{\nabla^1 q_h^n - \nabla^1 q_h^{n-1}}{c_2^n} = \frac{1}{c_2^n} \bigg(\frac{q_h^n - q_h^{n-1}}{c_1^n} - \frac{q_h^{n-1} - q_h^{n-2}}{c_1^{n-1}}\bigg),
    \end{align*}
    we have
    \begin{equation}\label{EQ:omega2-collected}
    \begin{aligned}
        &\nabla^2 q_h^n - \frac{1}{c_2^n} \left( \frac{c_2^n}{c_1^n} - \frac{c_2^{n-1}}{c_1^{n-1}} \right) \nabla^2 u_h^n + \Big( 3 - \frac{c_2^{n-1} c_3^n}{c_1^{n-1} c_2^n} \Big) \nabla^3 u_h^n 
        \\ &\qq= 
        - \frac{1}{c_2^n} \nabla^2 u_h^{n} \bigg( \frac{c_2^n}{c_1^n} - \frac{c_2^{n-1}}{c_1^{n-1}} + 1 \bigg)
        + \frac{1}{c_2^n} \nabla^2 u_h^{n-1} \bigg( \frac{c_1^{n-1}}{c_1^{n}} + 1 \bigg)
        - \frac{1}{c_2^n} \frac{c_1^{n-2}}{c_1^{n-1}} \nabla^2 u_h^{n-2} 
        \\ &\qqq
        + \bigg( 2 - \frac{c_2^{n-1} c_3^n}{c_1^{n-1} c_2^n}\bigg) \nabla^3 u_h^{n}
        + \frac{c_2^{n-1}}{c_1^{n}} \nabla^3 u_h^{n-1}
        - \frac{c_1^{n-2} c_2^{n-2}}{c_1^{n-1} c_2^n} \nabla^3 u_h^{n-2}
        \\ &\qq= 
        - \bigg( \frac{c_1^{n-2}}{c_2^n} + \frac{c_2^{n-1}}{c_1^n} \bigg) c_4^{n} \nabla^4 u_h^{n}
        + \frac{c_1^{n-2} c_2^{n-2} c_4^{n-1}}{c_1^{n-1} c_2^n} \nabla^4 u_h^{n-1}
        \\ &\qq= 
        \left(\frac{c_1^{n-2} c_2^{n-2} c_4^{n-1}}{c_1^{n-1} c_2^n} - \bigg( \frac{c_1^{n-2}}{c_2^n} + \frac{c_2^{n-1}}{c_1^n} \bigg) c_4^{n}\right)\nabla^4 u_h^{n}
        - \frac{c_1^{n-2} c_2^{n-2} c_4^{n-1} c_5^n}{c_1^{n-1} c_2^n} \nabla^5 u_h^{n}.
    \end{aligned}
    \end{equation}
    Collecting the terms with $\omega_1^{n}$ and resolving $\nabla^1 q_h^n$ by
    \begin{align*}
        \nabla^1 q_h^n = \frac{q_h^n - q_h^{n-1}}{c_1^n},
    \end{align*}
    we have
    \begin{equation}\label{EQ:omega1-collected}
    \begin{aligned}
        &\nabla^1 q_h^n + \Big( 2 - \frac{c_2^n}{c_1^n} \Big) \nabla^2 u_h^n + \big( 2 c_2^n - c_1^n \big) \nabla^3 u_h^n
        \\ &\qq=
        \frac{c_1^{n-1}}{c_1^n} \Big( \nabla^2 u_h^{n-1} - \nabla^2 u_h^n \Big)
        + c_1^{n-1} \nabla^3 u_h^n
        + \frac{c_1^{n-1} c_2^{n-1}}{c_1^n} \nabla^3 u_h^{n-1}
        \\ &\qq=
        - \frac{c_1^{n-1} c_2^{n-1}}{c_1^n} \, c_4^n \, \nabla^4 u_h^n 
    \end{aligned}
    \end{equation}
    In total, inserting \eqref{EQ:omega2-collected} and \eqref{EQ:omega1-collected} in \eqref{EQ:RR3-tmp1}, we conclude
    \begin{align*}
        r_h^{(3)}(t)
        &= \FF^{(2)}(t) - f(t)
        + A_4(t) \nabla^4 u_h^n
        + A_5(t) \nabla^5 u_h^n.
    \end{align*}
    where $A_4$ and $A_5$ are given in \eqref{EQ:Def-A4-A5}.
\end{proof}

For uniform time-steps, we have $c_m^n = m \tau$, and the coefficients in Lemma~\ref{LEM:residual-identities} simplify as follows.

\begin{Corollary}
    Let $r_h^{(k)}$ be the residual given in \eqref{EQ:def-residual}.
    Let $0 = t_0 < t_1 < t_2 < \dots < t_N = T$ be a uniform partition of $[0, T]$. 
    Then, it holds
    \begin{align}
        r_h^{(1)}(t) &= \FF^{(0)}(t) - f(t), 
        \\
        r_h^{(2)}(t) 
        &= \FF^{(1)}(t) - f(t) - 3 \tau\omega_1^{n}(t) \nabla^3 u_h^n,
        \\
        r_h^{(3)}(t) 
        &= 
        \FF^{(2)}(t) - f(t) - \big( 6 \tau \omega_2^{n}(t) + 8 \tau^2 \omega_1^{n}(t) \big) \nabla^4 u_h^n
        - 20 \tau^2 \omega_2^{n}(t) \nabla^5 u_h^n
    \end{align}
    for $t \in I_n$, $n \ge k$, where the values $u_h^{-\ell}$ at negative times $t_{-\ell}$ for $\ell = 1, \dots, k-1$ are defined by \eqref{EQ:def-neg-time-vals-u}.
\end{Corollary}

For \BDF{4}, we restrict the analysis of the residual to the uniform time-step case in the following lemma.

\begin{Lemma}\label{LEM:BDF4-residual}
    Let $r_h^{(4)}$ be the residual given in \eqref{EQ:def-residual} for $k=4$.
    Let $0 = t_0 < t_1 < t_2 < \dots < t_N = T$ be a uniform partition of $[0, T]$. 
    Then, it holds
    \begin{align*}
        r_h^{(4)}(t) &= \FF^{(3)}(t) - f(t)
        - 10 \tau  \Big[3 \tau^2 \omega_1^{n}(t) + 2 \tau \omega_2^{n}(t) + \omega_3^{n}(t) \Big]\nabla^5u_h^n
        \\&\qq
        - 10 \tau^2 \Big[9 \tau \omega_2^{n}(t) + 4  \omega_3^{n}(t) \Big] \nabla^6u_h^n
        - 210\tau^3\,\omega_3^{n}(t)\,\nabla^7u_h^n
    \end{align*}
    for $t \in I_n$, $n \ge 4$, where the values $u_h^{-\ell}$ at negative times $t_{-\ell}$ for $\ell = 1, \dots, k-1$ are defined by \eqref{EQ:def-neg-time-vals-u}.
\end{Lemma}
\begin{proof}
    For a uniform mesh, we have $c_m^{n-j} = m\tau$ for every shift $0 \le j \le n$ and every $m \ge 1$, where we use the negative times $t_{-\ell} = -\ell \tau$ if $m > n-j$.
    In particular, using \eqref{EQ:Def-vhnk} and \eqref{EQ:def-qhn} we obtain
    \begin{align}\label{EQ:def-qhn-BDF4}
        q_h^{n} = -\tau\nabla^2u_h^{n} - 2\tau^2\nabla^3u_h^{n} - 6\tau^3\nabla^4u_h^{n}.
    \end{align}
    For a uniform mesh, we obtain for the nested divided differences
    \begin{align}
        \label{EQ:binomial-identity}
        \nabla^i\big(\nabla^ju_h^n\big) = \binom{i+j}{i}\nabla^{i+j}u_h^n.
    \end{align}
    In particular $\nabla^1(\nabla^1u_h^n)=2\nabla^2u_h^n$, $\nabla^2(\nabla^1u_h^n)=3\nabla^3u_h^n$, $\nabla^3(\nabla^1u_h^n)=4\nabla^4u_h^n$.
    Applying \eqref{EQ:binomial-identity} to each term of $q_h^{n}$ in \eqref{EQ:def-qhn-BDF4} gives
    \begin{align*}
        \nabla^1q_h^n &= 
        -3\tau\nabla^3u_h^n - 8\tau^2\nabla^4u_h^n - 30\tau^3\nabla^5u_h^n, \\
        \nabla^2q_h^n &= 
        -6\tau\nabla^4u_h^n - 20\tau^2\nabla^5u_h^n - 90\tau^3\nabla^6u_h^n, \\
        \nabla^3q_h^n &= 
        -10\tau\nabla^5u_h^n - 40\tau^2\nabla^6u_h^n - 210\tau^3\nabla^7u_h^n.
    \end{align*}
    Using $\omega_1^{n}(t)^2=\omega_2^{n}(t)-\tau\omega_1^{n}(t)$ and $\omega_1^{n}(t)\omega_2^{n}(t)=\omega_3^{n}(t)-2\tau\omega_2^{n}(t)$ to eliminate powers of $\omega_1^{n}(t)$ from Lemma~\ref{LEM:rho_i}, we rewrite $\rho_2^n,\rho_3^n,\rho_4^n$ by
    \begin{align}
        \rho_2^n(t) &= 2\,\omega_1^{n}(t), \\
        \rho_3^n(t) &= 3\,\omega_2^{n}(t) + 3\tau\,\omega_1^{n}(t), \\
        \rho_4^n(t) &= 4\,\omega_3^{n}(t) + 6\tau\,\omega_2^{n}(t) + 8\tau^2\,\omega_1^{n}(t).
    \end{align}
    Substituting all of the above into \eqref{EQ-RR-master} for $k=4$ and collecting by $\omega_1^{n}(t),\omega_2^{n}(t),\omega_3^{n}(t)$ yields
    \begin{align}
        r_h^{(4)}(t) = \FF^{(3)}(t) - f(t) + \omega_1^{n}(t)\,X_1 + \omega_2^{n}(t)\,X_2 + \omega_3^{n}(t)\,X_3,
    \end{align}
    where we obtain by canceling the lower-order divided differences with \eqref{EQ:binomial-identity}
    \begin{align}
        X_1 &\coloneqq \nabla^1q_h^n - \nabla^1(\nabla^1u_h^n) + 2\nabla^2u_h^n + 3\tau\nabla^3u_h^n + 8\tau^2\nabla^4u_h^n = -30\tau^3\,\nabla^5u_h^n, \\
        X_2 &\coloneqq \nabla^2q_h^n - \nabla^2(\nabla^1u_h^n) + 3\nabla^3u_h^n + 6\tau\nabla^4u_h^n = -20\tau^2\,\nabla^5u_h^n - 90\tau^3\,\nabla^6u_h^n, \\
        X_3 &\coloneqq \nabla^3q_h^n - \nabla^3(\nabla^1u_h^n) + 4\nabla^4u_h^n = -10\tau\,\nabla^5u_h^n - 40\tau^2\,\nabla^6u_h^n - 210\tau^3\,\nabla^7u_h^n.
    \end{align}
\end{proof}

Estimating the identities in Lemma~\ref{LEM:residual-identities} and Lemma~\ref{LEM:BDF4-residual}, we obtain the following lemma.

\begin{Lemma}\label{LEM:est-residual-final}
    Let $r_h^{(k)}$ be the residual given in \eqref{EQ:def-residual} for $k \in \{1, 2, 3, 4\}$ and $\star\in\{H^{-1},L^2\}$.
    Then, it holds
    \begin{align*}
        \int_J\norm{r_h^{(k)}}_\star^2\dt \le k\Big(\FF^{(k)}(\star)+\NN^{(k)}(\star)\Big),
    \end{align*}
    where $\FF^{(k)}(\star)$ is defined in \eqref{EQ:def-data-approx-err-est} and $\NN^{(k)}(\star)$ is defined in \eqref{EQ:def-NN1-NN2-NN3} for \VSBDF{k}, $k\in \{1,2,3\}$, and in \eqref{EQ:def-NN4-uni} for uniform \BDF{4}.
\end{Lemma}
\begin{proof}
    Fix $n\ge k$ and $t\in I_n$. By Lemma~\ref{LEM:residual-identities} (for $k\le3$) and Lemma~\ref{LEM:BDF4-residual} (for $k=4$, uniform mesh), write
    \begin{align*}
        r_h^{(k)}(t) = T_0(t) + T_1(t) + \dots + T_{k-1}(t),
        \qquad T_0(t)\coloneqq\FF^{(k-1)}(t)-f(t),
    \end{align*}
    with $T_1,\dots,T_{k-1}$ the individual divided-difference terms.
    Then, by the discrete Cauchy--Schwarz inequality we obtain
    \begin{align}\label{EQ:discr-CS}
        \norm{r_h^{(k)}(t)}_\star^2 \le k\sum_{i=0}^{k-1}\norm{T_i(t)}_\star^2.
    \end{align}
    Integrating over $I_n$ and summing $n=k,\dots,N$ gives
    \begin{align*}
        \int_J\norm{r_h^{(k)}}_\star^2\dt
        \le k\int_J\norm{T_0}_\star^2\dt + k\sum_{i=1}^{k-1}\int_J\norm{T_i}_\star^2\dt
        = k\,\FF^{(k)}(\star) + k\,\NN^{(k)}(\star).
    \end{align*}
    It remains to evaluate each $\NN^{(k)}(\star)$ for each $k$.

    \textbf{Case $k=1$:} The terms vanish, yielding
    \begin{align}\label{EQ:EE1}
        \NN^{(1)}(\star)=0.
    \end{align}

    \textbf{Case $k=2$:} 
    We have $T_1(t) = -\dfrac{c_1^{n-1}c_3^n}{c_1^n}\,\omega_1^{n}(t)\,\nabla^3u_h^n$. 
    Using \eqref{EQ:mom-11} yields
    \begin{align*}
        \int_{I_n}\norm{T_1(t)}_\star^2\dt
        = \frac{\tau_n^3}{3}\Big(\frac{c_1^{n-1}c_3^n}{c_1^n}\Big)^2\norm{\nabla^3u_h^n}_\star^2
        = \frac{(c_1^{n-1})^2(c_3^n)^2}{3}\,c_1^n\,\norm{\nabla^3u_h^n}_\star^2.
    \end{align*}
    Hence
    \begin{align}
        \label{EQ:E2}
        \NN^{(2)}(\star) = \sum_{n=2}^{N} \frac{(c_1^{n-1})^2(c_3^n)^2}{3}\,c_1^n\,
        \norm{\nabla^3u_h^n}_\star^2.
    \end{align}

    \textbf{Case $k=3$:} 
    We have $T_1(t)=A_4(t)\nabla^4u_h^n$ and $T_2(t)=A_5(t)\nabla^5u_h^n$.
    For convenience, we abbreviate $A_4(t)=\alpha_n\omega_1^{n}(t)+\beta_n\omega_2^{n}(t)$ and $A_5(t)=\gamma_n\omega_2^{n}(t)$,~cf.~\eqref{EQ:Def-A4-A5},
    \begin{align*}
        \alpha_n &= -\frac{c_1^{n-1}c_2^{n-1}c_4^n}{c_1^n}, &
        \beta_n &= \frac{c_1^{n-2}c_2^{n-2}c_4^{n-1}}{c_1^{n-1}c_2^n}
            - \Big(\frac{c_1^{n-2}}{c_2^n}+\frac{c_2^{n-1}}{c_1^n}\Big)c_4^n, &
        \gamma_n &= -\frac{c_1^{n-2}c_2^{n-2}c_4^{n-1}c_5^n}{c_1^{n-1}c_2^n}.
    \end{align*}
    By \eqref{EQ:mom-11}--\eqref{EQ:mom-22}, we obtain
    \begin{align*}
        \int_{I_n}A_4(t)^2\dt = \frac{\tau_n^3}{3}\alpha_n^2 + \frac{\tau_n^4}{6}\alpha_n\beta_n
        + \frac{\tau_n^5}{30}\beta_n^2,
        \qquad
        \int_{I_n}A_5(t)^2\dt = \frac{\tau_n^5}{30}\gamma_n^2.
    \end{align*}
    Thus, we conclude 
    \begin{align}
        \label{EQ:E3}
        \NN^{(3)}(\star) = \sum_{n=3}^{N}\Big[\Big(\frac{\tau_n^3}{3}\alpha_n^2
        + \frac{\tau_n^4}{6}\alpha_n\beta_n + \frac{\tau_n^5}{30}\beta_n^2\Big)
        \norm{\nabla^4u_h^n}_\star^2
        + \frac{\tau_n^5}{30}\gamma_n^2\,\norm{\nabla^5u_h^n}_\star^2\Big].
    \end{align}

    \textbf{Case $k=4$ (uniform time-steps):} 
    We have
    \begin{align*}
        T_1(t) &= -10\tau\big[3\tau^2\omega_1^{n}(t)+2\tau\omega_2^{n}(t)+\omega_3^{n}(t)\big]\nabla^5u_h^n,
        \\
        T_2(t) &= -10\tau^2\big[9\tau\omega_2^{n}(t)+4\omega_3^{n}(t)\big]\nabla^6u_h^n,
        \\
        T_3(t) &= -210\tau^3\,\omega_3^{n}(t)\,\nabla^7u_h^n.
    \end{align*}
    Integrating each bracket's square yields
    \begin{align*}
        \int_{I_n}\big[3\tau^2\omega_1^{n}+2\tau\omega_2^{n}+\omega_3^{n}\big]^2\dt
        &= \frac{1073}{210}\,\tau^7,
        \\
        \int_{I_n}\big[9\tau\omega_2^{n}+4\omega_3^{n}\big]^2\dt
        &= \frac{1579}{210}\,\tau^7,
        \\
        \int_{I_n}(\omega_3^{n})^2\dt &= \frac{8\tau^7}{105}.
    \end{align*}
    Multiplying by the outer prefactors, we conclude
    \begin{align}
        \label{EQ:E4}
        \NN^{(4)}(\star) = \sum_{n=4}^N\Big[\frac{10730}{21}\tau^9\norm{\nabla^5u_h^n}_\star^2
        + \frac{15790}{21}\tau^{11}\norm{\nabla^6u_h^n}_\star^2
        + 3360\,\tau^{13}\norm{\nabla^7u_h^n}_\star^2\Big].
    \end{align}
\end{proof}

The constant $k$ in Lemma~\ref{LEM:est-residual-final} results from a rather crude estimate, which can be improved as follows.

\begin{Remark}\label{REM:factor-k}
    The factor $k$ in Lemma~\ref{LEM:est-residual-final} comes from the discrete Cauchy--Schwarz inequality applied to the $k$ terms $T_0, \dots, T_{k-1}$ in \eqref{EQ:discr-CS}.
    However, this is not optimal for $k\ge3$ since, for sufficiently smooth solutions, $T_2, \dots, T_{k-1}$ are of higher order. 
    Instead, applying Young's inequality with a parameter $\varepsilon > 0$ yields
    \begin{align*}
        \norm{r_h^{(k)}(t)}_\star^2 \le 
        2 (1 + \varepsilon) \big(\norm{T_0(t)}_\star^2 + \norm{T_1(t)}_\star^2\big)
        + (1 + \varepsilon^{-1}) (k-2) \sum_{i=2}^{k-1}\norm{T_i(t)}_\star^2,
    \end{align*}
    i.e.,~the constant for the asymptotically dominant term is $2 (1 + \varepsilon)$ instead of $k$.
    For simplicity, we keep the factor $k$. 
\end{Remark}

\subsection{Estimating $\delta^{(k)}$}
Having estimated the residual, we now bound the elliptic-reconstruction
defect $\delta^{(k)}$ given in \eqref{EQ:ell-rec-defect} and its time derivative, which yields the estimators $\Theta^{(k)}$, $\Theta_\infty^{(k)}$, and $\Lambda^{(k)}$ from Definition~\ref{DEF:collecting-estimators}.

We start by resolving $\nabla^i u_h^n$ for $i \ge 0$ into $\nabla^0 u_h^{n-j} = u_h^{n-j}$, via coefficients $\delta_j^{i,n}$ defined by
\begin{align}\label{EQ:def-delta-coeff}
    \nabla^i u_h^n = \sum_{j=0}^{i} \delta_j^{i, n} \, u_h^{n-j}.
\end{align}
Substituting the recursive definition of $\nabla^i$ \eqref{DEF:backward-div-diff} and matching coefficients gives the recursion
\begin{align}\label{EQ:def-delta-coeff-recursion}
    \delta_j^{i, n} = \frac{1}{c_i^n}\Big[\delta_j^{i-1, n} - \delta_{j-1}^{i-1, n-1}\Big],
    \qquad j = 0, \ldots, i,
\end{align}
with the convention $\delta_j^{i-1, n} \coloneqq 0$ for $j < 0$ or $j > i-1$, and base case $\delta_0^{0, n} = 1$.
With these coefficients, we obtain the following pointwise bound for $\delta^{(k)}$.

\begin{Lemma}\label{LEM:delta-est}
    Let $\star \in \{L^2, H^1\}$.
    For $\delta^{(k)} = \recW{k} - \recM{k}$, it holds for $t \in I_n$, $n \ge k$,
    \begin{align*}
        \norm{\delta^{(k)}(t)}_\star^2
        &\le 
        (k+1) \sum_{j=0}^{k} \big(A_j^{n,(k)}\big)^2\, \eta(u_h^{n-j}; \star),
    \end{align*}
    where the coefficients are defined by
    \begin{align}\label{EQ:alpha-coeff}
        A_j^{n,(k)} \coloneqq \sup_{t \in I_n} \big|\alpha_j^{n,(k)}(t)\big|,
        \qquad
        \alpha_j^{n,(k)}(t) \coloneqq \sum_{i=j}^{k} \omega_i^{n}(t)\, \delta_j^{i, n}, \qq n \ge k,
    \end{align}
    and $A_j^{n,(k)}\coloneqq 0$ if $n\notin\{k,\dots,N\}$
\end{Lemma}
\begin{proof}
    Using \eqref{EQ:def-delta-coeff}, we have
    \begin{align*}
        \delta^{(k)}(t)
        =
        \sum_{i=0}^{k} \omega_i^{n}(t) (\RR - I) \nabla^i u_h^{n}
        =
        \sum_{i=0}^{k} \omega_i^{n}(t) \sum_{j=0}^{i} \delta_j^{i, n} \, (\RR - I) u_h^{n-j}
    \end{align*}
    Next, define
    \begin{align*}
        \alpha_j^{n,(k)}(t) \coloneqq \sum_{i=j}^{k} \omega_i^{n}(t)\, \delta_j^{i, n}.
    \end{align*}
    Exchanging the order of summation, we obtain
    \begin{align*}
        \delta^{(k)}(t)
        &=
        \sum_{j=0}^{k} \left[\sum_{i=j}^{k} \omega_i^{n}(t)\, \delta_j^{i, n}\right] (\RR - I) u_h^{n-j}
        =
        \sum_{j=0}^{k} \alpha_j^{n,(k)}(t)\, (\RR - I) u_h^{n-j}.
    \end{align*}
    By the discrete Cauchy--Schwarz inequality applied to the sum of $k+1$ terms,
    \begin{align*}
        \norm{\delta^{(k)}(t)}_\star^2
        \le
        (k+1) \sum_{j=0}^{k} \big|\alpha_j^{n,(k)}(t)\big|^2\, \norm{(\RR - I) u_h^{n-j}}_\star^2, \qquad t \in I_n.
    \end{align*}
    Applying the \emph{a~posteriori} bound \eqref{EQ:ell-rec-apost-L2} or \eqref{EQ:ell-rec-apost-H1}, respectively, and bounding $|\alpha_j^{n,(k)}(t)| \le A_j^{n,(k)}$, we conclude the proof.
\end{proof}

Similar to \eqref{EQ:def-delta-coeff}, we continue by rewriting $\nabla^i u_h^n$ for $i \ge 1$ into $\nabla^1 u_h^{n-j}$ terms via coefficients $\gamma_j^{i,n}$ defined by
\begin{align}\label{EQ:def-gamma-coeff}
    \nabla^i u_h^n = \sum_{j=0}^{i-1} \gamma_j^{i, n} \, \nabla^1 u_h^{n-j}.
\end{align}
Substituting the recursive definition of $\nabla^i$ \eqref{DEF:backward-div-diff} and matching coefficients gives the recursion
\begin{align*}
    \gamma_j^{i, n} = \frac{1}{c_i^n}\Big[\gamma_j^{i-1, n} - \gamma_{j-1}^{i-1, n-1}\Big],
    \qquad j = 0, \ldots, i-1,
\end{align*}
with the convention $\gamma_j^{i-1, n} \coloneqq 0$ for $j < 0$ or $j > i-2$, and $\gamma_0^{1, n} = 1$.
In particular, for $i=2$ we have $\gamma_0^{2, n} = 1/c_2^n$, $\gamma_1^{2, n} = -1/c_2^n$.
These coefficients allow us to bound $\ptl_t\delta^{(k)}$ analogously to $\delta^{(k)}$ itself, resolving it into $(\RR-I)\nabla^1u_h^{n-j}$ instead of $(\RR-I)u_h^{n-j}$.

\begin{Lemma}\label{LEM:dt-delta-est}
    For $\delta^{(k)} = \recW{k} - \recM{k}$ and $\star \in \{H^{-1}, L^2\}$, it holds
    \begin{align*}
        \int_{t_{k-1}}^{t_N} \norm{\ptl_t \delta^{(k)}}_{\star}^2 \dt
        &\le 
        k \sum_{m=1}^{N} \Xi_m^{(k)}\, \eta(\nabla^1 u_h^{m}; \star)
        \eqqcolon
        \Lambda^{(k)}(\star),
    \end{align*}
    where the coefficients are defined by 
    \begin{align}\label{EQ:coeff-beta}
        \Xi_m^{(k)} \coloneqq \sum_{j=0}^{k-1} \tau_{m+j}\, \big(B_j^{m+j,(k)}\big)^2,
        \qq
        B_j^{n, (k)} \coloneqq \sup_{t \in I_n} |\beta_j^{n, (k)}(t)|,
        \qq
        \beta_j^{n, (k)}(t) = 
        \sum_{i=j+1}^{k} \ptl_t \omega_i^n(t)\, \gamma_j^{i, n}.
    \end{align}
\end{Lemma}
\begin{proof}
    By definition,
    \begin{align*}
        \ptl_t \delta^{(k)}(t)
        =
        \sum_{i=1}^{k}
        \ptl_t \omega_i^n(t)\,
        (\RR - I)
        \nabla^i u_h^{n}, \qquad t \in I_n.
    \end{align*}
    Using \eqref{EQ:def-gamma-coeff}
    and exchanging the order of summation, we obtain
    \begin{align*}
        \ptl_t \delta^{(k)}(t)
        &=
        \sum_{i=1}^{k} \ptl_t \omega_i^n(t) (\RR-I) \sum_{j=0}^{i-1} \gamma_j^{i, n}\, \nabla^1 u_h^{n-j}
        =
        \sum_{j=0}^{k-1}
        \beta_j^{n, (k)}(t) 
        (\RR - I) \nabla^1 u_h^{n-j}.
    \end{align*}
    Next, we have
    \begin{align*}
        \norm{\ptl_t \delta^{(k)}(t)}_\star^2
        \le
        k \sum_{j=0}^{k-1} \big|\beta_j^{n, (k)}(t)\big|^2\, \norm{(\RR - I)\nabla^1 u_h^{n-j}}_\star^2 
    \end{align*}
    for $t \in I_n$.
    With $B_j^{n, (k)} = \sup_{t \in I_n} |\beta_j^{n, (k)}(t)|$, integrating over $I_n$ and applying the \emph{a~posteriori} bound \eqref{EQ:ell-rec-apost-dual} or \eqref{EQ:ell-rec-apost-L2}, respectively, yields
    \begin{align*}
        \int_{I_n}\norm{\ptl_t \delta^{(k)}(t)}_\star^2 \dt
        \le 
        k \, \tau_n \sum_{j=0}^{k-1} \big|B_j^{n, (k)}\big|^2\, \eta(\nabla^1 u_h^{n-j}; \star).
    \end{align*}
    Summing over $n = k,\dots,N$, we obtain
    \begin{align*}
        \int_{t_{k-1}}^{t_N} \norm{\ptl_t \delta^{(k)}}_\star^2 \dt
        &\le 
        k \sum_{n=k}^{N} \tau_n \sum_{j=0}^{k-1} \big|B_j^{n,(k)}\big|^2\, \eta(\nabla^1 u_h^{n-j}; \star).
    \end{align*}
    Re-indexing $m \coloneqq n-j$ for each fixed $j$ (so $n=m+j$) and exchanging the order of summation and setting $B_j^{n,(k)} = 0$ if $n \notin \{k,\ldots,N\}$, this gives
    \begin{align*}
        \int_{t_{k-1}}^{t_N} \norm{\ptl_t \delta^{(k)}}_\star^2 \dt
        \le 
        k \sum_{m=1}^{N} \Xi_m^{(k)}\, \eta(\nabla^1 u_h^{m}; \star).
    \end{align*}
\end{proof}

For uniform time-steps, the weights $\Xi_m^{(k)}$ in Lemma~\ref{LEM:dt-delta-est} can be computed explicitly.

\begin{Remark}
    For uniform time-steps, we have
    \begin{align*}
        \Xi_m^{(2)} = \frac{5}{2} \tau,
        \qq
        \Xi_m^{(3)} = \frac{29}{6} \tau
    \end{align*}
    for $k \le m \le N - k + 1$.
\end{Remark}

Moreover, the coefficients $\beta_j^{n,(k)}(t_n)$ are closely related to the discrete convolution kernels used in the stability analysis of variable-step BDF methods.

\begin{Remark}
    The coefficients $\beta_j^{n, (k)}$ defined in \eqref{EQ:coeff-beta} evaluated at $t_n$ correspond precisely to the discrete convolution kernels of the variable-step \BDF{k} formula 
    introduced in~\cite{LiLiao:2022}, namely,
    \begin{align*}
        \beta_j^{n,(k)}(t_n) = d_j^{(k,n)},\qq n \ge k, \quad  k - 1 \ge j \ge 0,
    \end{align*}
    where $d_j^{(k,n)}$ are given explicitly in \cite[Eqs.~(1.9)--(1.10)]{LiLiao:2022} for $k = 2, 3$.
\end{Remark}

We now continue with estimating $\delta^{(k)}$, deriving its $L^2$-in-time bound $\Theta^{(k)}$.

\begin{Lemma}\label{LEM:grad-delta-est}
    Define, for $\star\in\{L^2,H^1\}$,
    \begin{align*}
        \Theta^{(k)}(\star) \coloneqq  (k+1)\sum_{m=0}^{N} \Psi_m^{(k)}\,\eta(u_h^m; \star),
        \qquad
        \Psi_m^{(k)} \coloneqq \sum_{j=0}^{k} \tau_{m+j}\,\big(A_j^{m+j,(k)}\big)^2,
    \end{align*}
    where the coefficients $A_j^{n,(k)}$ are defined in \eqref{EQ:alpha-coeff}. Then, it holds
    \begin{align*}
        \int_{t_{k-1}}^{t_N} \norm{\delta^{(k)}}_\star^2 \dt \le \Theta^{(k)}(\star).
    \end{align*}
\end{Lemma}
\begin{proof}
    By Lemma~\ref{LEM:delta-est}, we have for $t\in I_n$, $n\ge k$,
    \begin{align*}
        \norm{\delta^{(k)}(t)}_\star^2 \le  (k+1)\sum_{j=0}^{k}\big(A_j^{n,(k)}\big)^2\,
        \eta(u_h^{n-j};\star).
    \end{align*}
    Since $A_j^{n,(k)}$ is independent of $t$, integrating over $I_n$ yields
    \begin{align*}
        \int_{I_n}\norm{\delta^{(k)}(t)}_\star^2\dt
        \le (k+1)\,\tau_n\sum_{j=0}^{k}\big(A_j^{n,(k)}\big)^2\,\eta(u_h^{n-j};\star).
    \end{align*}
    Summing over $n=k,\dots,N$ and changing the order of summation yields
    \begin{align*}
        \int_{t_{k-1}}^{t_N}\norm{\delta^{(k)}}_\star^2\dt
        &\le (k+1)\sum_{n=k}^{N}\tau_n\sum_{j=0}^{k}\big(A_j^{n,(k)}\big)^2\,\eta(u_h^{n-j};\star)
        \\
        &= (k+1)\sum_{m=0}^N\Psi_m^{(k)}\,\eta(u_h^m;\star).
    \end{align*}
\end{proof}

Similarly, we obtain the corresponding $L^\infty$ in-time bound.

\begin{Lemma}\label{LEM:delta-Linfty-est}
    Define, for $\star\in\{L^2,H^1\}$,
    \begin{align*}
        \Theta_\infty^{(k)}(\star) \coloneqq (k+1)\max_{k\le n\le N}\sum_{j=0}^{k}
        \big(A_j^{n,(k)}\big)^2\,\eta(u_h^{n-j};\star),
    \end{align*}
    where the coefficients $A_j^{n,(k)}$ are defined in \eqref{EQ:alpha-coeff}.
    Then, it holds
    \begin{align*}
        \norm{\delta^{(k)}}_{L^\infty(t_{k-1}, t_N; \star)}^2 \le \Theta_\infty^{(k)}(\star).
    \end{align*}
\end{Lemma}
\begin{proof}
    By Lemma~\ref{LEM:delta-est}, for every $t\in I_n$, $n\ge k$,
    \begin{align*}
        \norm{\delta^{(k)}(t)}_\star^2 \le (k+1)\sum_{j=0}^{k}\big(A_j^{n,(k)}\big)^2\,\eta(u_h^{n-j}; \star).
    \end{align*}
    Since the right-hand side does not depend on $t$, taking the supremum over $t\in I_n$ yields
    \begin{align*}
        \sup_{t\in I_n}\norm{\delta^{(k)}(t)}_\star^2 \le (k+1)\sum_{j=0}^{k}\big(A_j^{n,(k)}\big)^2\,
        \eta(u_h^{n-j}; \star).
    \end{align*}
    Thus, we obtain for the supremum over $(t_{k-1}, t_N)$
    \begin{align*}
        \norm{\delta^{(k)}}_{L^\infty(t_{k-1}, t_N;\star)}^2
        = \max_{k\le n\le N}\ \sup_{t\in I_n}\norm{\delta^{(k)}(t)}_\star^2
        \le (k+1)\max_{k\le n\le N}\sum_{j=0}^{k}\big(A_j^{n,(k)}\big)^2\,\eta(u_h^{n-j};\star)
        = \Theta_\infty^{(k)}(\star).
    \end{align*}
\end{proof}

We can also conclude a less sharp but decoupled estimate for $\Theta_\infty^{(k)}(\star)$.

\begin{Corollary}\label{COR:decoupled-Linfty}
    Define the maximum constants
    \begin{align*}
        \widehat A_j^{(k)} \coloneqq \max_{k\le n\le N} A_j^{n,(k)}, 
    \end{align*}
    for $j=0,\dots,k$, where the coefficients $A_j^{n,(k)}$ are defined in \eqref{EQ:alpha-coeff}. 
    Then
    \begin{align*}
        \Theta_\infty^{(k)}(\star) \le (k+1)\Big(\sum_{j=0}^k\big(\widehat A_j^{(k)}\big)^2\Big)\,
        \max_{0\le m\le N}\eta(u_h^m;\star).
    \end{align*}
\end{Corollary}
\begin{proof}
    Shifting the max inside the sum yields
    \begin{align*}
        \Theta_\infty^{(k)}(\star) \le (k+1)\sum_{j=0}^k \max_{k\le n\le N}\Big[\big(A_j^{n,(k)}\big)^2\,
        \eta(u_h^{n-j};\star)\Big].
    \end{align*}
    For each fixed $j$, we have
    \begin{align*}
        \max_{k\le n\le N}\Big[\big(A_j^{n,(k)}\big)^2\eta(u_h^{n-j};\star)\Big]
        \le 
        \big(\widehat A_j^{(k)}\big)^2\max_{k\le n\le N}\eta(u_h^{n-j};\star).
    \end{align*}
    Finally, we conclude 
    \begin{align*}
        \max_{k\le n\le N}\eta(u_h^{n-j};\star) 
        \le \max_{0\le m\le N}\eta(u_h^m;\star).
    \end{align*}
    Summing over $j=0,\dots,k$ concludes the proof.
\end{proof}

\subsection{Estimating $\rho^{(k)}$ and $\rho_h^{(k)}$}
It remains to bound the inter-order terms $\rho^{(k)}$ and $\rho_h^{(k)}$,
which control the difference between the order $(k-1)$ and order $k$ reconstructions, and yield the estimator $\Sest^{(k)}$.

\begin{Lemma}\label{LEM:est-rho}
    For $\rho_h^{(k)} = \recM{k-1} - \recM{k}$, it holds
    \begin{align}
        \int_{t_{k-1}}^{t_N} \norm{\Delta_h \rho_h^{(k)}}^2 \dt &\le \sum\limits_{n = k}^N \omega^{n, (k)} \norm{\nabla^k (\Delta_hu_h^n)}^2 = \Sest^{(k)}(H^2)
        \label{EQ:rho-h-est-2},
    \end{align}
    where $\omega^{n, (k)} = \int_{t_{n-1}}^{t_n} |\omega_k^{n}(t)|^2 \dt$.
    Moreover, it holds for $\rho^{(k)} = \recW{k-1} - \recW{k}$
    \begin{align}
        \int_{t_{k-1}}^{t_N} \norm{\grad \rho^{(k)}}^2 \dt &\le \sum\limits_{n = k}^N \omega^{n, (k)} \Big(  \norm{\nabla^k (\grad u_h^n)}^2 + \eta(\nabla^k u_h^n; H^1) \Big) = \Sest^{(k)}(H^1) \label{EQ:rho-est-1}.
    \end{align}
\end{Lemma}
\begin{proof}
    We start with \eqref{EQ:rho-h-est-2}. We have by \eqref{EQ:Newton-interpol-poly} and \eqref{EQ:k+1-time-rec-def} 
    \begin{align*}
        \rho_h^{(k)} = \recM{k-1} - \recM{k} = - \omega_k^{n}(t) \nabla^k u_h^n
    \end{align*}
    for $t \in I_n$.
    Thus 
    \begin{align*}
        \int_{t_{n-1}}^{t_n} \norm{\Delta_h \rho_h^{(k)}}^2 \dt = \omega^{n, (k)} \norm{\nabla^k (\Delta_hu_h^n)}^2.
    \end{align*}
    Summing over $n=k,\dots,N$ gives \eqref{EQ:rho-h-est-2} with $\Sest^{(k)}(H^2)$.

    For \eqref{EQ:rho-est-1}, recall that $\rho^{(k)} = \RR \rho_h^{(k)}$ with $\rho_h^{(k)}(t) \in V_h$.
    For any $v_h \in V_h$, the definitions \eqref{EQ: discrete Laplacian def} and \eqref{EQ:ell-reconstruction} yield
    \begin{align*}
        a\big((\RR - I) v_h, v_h\big) = \dual{-\Delta_h v_h, v_h} - a(v_h, v_h) = 0,
    \end{align*}
    and thus $\norm{\grad \RR v_h}^2 = \norm{\grad v_h}^2 + \norm{\grad (\RR - I) v_h}^2$.
    Applying this and using \eqref{EQ:ell-rec-apost-H1}, we obtain
    \begin{align*}
        \norm{\grad\rho^{(k)}}^2 \le \norm{\grad\rho_h^{(k)}}^2 + \eta(\rho_h^{(k)}; H^1).
    \end{align*}
    Since $\rho_h^{(k)}(t) = -\omega_k^n(t)\nabla^k u_h^n$ for $t \in I_n$, integrating over $I_n$ yields
    \begin{align*}
        \int_{t_{n-1}}^{t_n} \norm{\grad\rho^{(k)}}^2 \dt \le \omega^{n, (k)} \Big( \norm{\nabla^k (\grad u_h^n)}^2 + \eta(\nabla^k u_h^n; H^1) \Big).
    \end{align*}
    Finally, summing over $n=k,\dots,N$ gives \eqref{EQ:rho-est-1}.
\end{proof}

This completes the proofs of Corollary~\ref{COR:first-apost-Linfty} and Corollary~\ref{COR:second-apost-Linfty}.

\section{Numerical experiments}\label{SEC:numerics}

In this section, we compute the error estimators on $\Omega=(0,1)$ with homogeneous Dirichlet boundary conditions, discretized in space by continuous piecewise linear ($\mathbb P_1$) or quadratic ($\mathbb P_2$) finite elements on uniform partitions with mesh size $h$, and in time by \BDF{k}, $k=1,\dots,4$ with uniform step-size $\tau$, as well as by \VSBDF{k}, $k = 1, 2, 3$, with alternating step-sizes.

\subsection{Setup}\label{SUBSEC:setup}

We consider two examples: Example~1 for which the temporal error dominates, and Example~2 where the spatial error dominates.

\newcounter{numexample}
\subsubsection{Example~1.}\refstepcounter{numexample}\label{SEC:ExampleTime}
As a first example, we consider
\begin{align}\label{EQ:ExampleTime}
    u(x,t) = \ee^{- (t-t_0)^2/(2 \sigma^2)}\sin(\pi x), \qquad (x,t)\in\Omega\times(0,T],
\end{align}
with $t_0 = 0.5$, $\sigma = 0.05$, $T=1$. Since $-\Delta u = \pi^2 u$ and $\partial_t u = -\tfrac{t-t_0}{\sigma^2}u$, the
right-hand side is given by
\begin{align}\label{EQ:ExampleTime-rhs}
    f(x,t) = (\pi^2-\frac{t - t_0}{\sigma^2})u(x,t),
\end{align}
and the initial data is $u_0(x) = u(x, 0)$. 

\subsubsection{Example~2.}\refstepcounter{numexample}\label{SEC:ExampleSpace}
As a second example, we consider the exact solution
\begin{align}\label{EQ:ExampleSpace}
    u(x,t) = \ee^{-t}\sin(\pi x), \qquad (x,t)\in\Omega\times(0,T],
\end{align}
with $T=1$. Since $-\Delta u = \pi^2 u$ and $\partial_t u = -u$, the
right-hand side is given by
\begin{align}\label{EQ:ExampleSpace-rhs}
    f(x,t) = (\pi^2-1)\, \ee^{-t}\sin(\pi x),
\end{align}
and the initial data is $u_0(x) = \sin(\pi x)$. 

\def\errI{\EE_I}
\def\errII{\EE_{II}}
\def\estI{\mathrm{Est}_I}
\def\estII{\mathrm{Est}_{II}}
\def\effI{\mathrm{Eff}_{I}}
\def\effII{\mathrm{Eff}_{II}}

\subsubsection{Implementation details}\label{SEC:impl-details}

For the numerical experiment, we compute the left-hand and right-hand sides of Corollary~\ref{COR:first-apost-Linfty} and Corollary~\ref{COR:second-apost-Linfty} and abbreviate them by
\begin{equation}\label{EQ:def-err-est}
\begin{aligned}
    \errI^2 &\coloneqq \|\varepsilon^{(k)}\|^2_{L^\infty(J;L^2(\Omega))}
        + \|\grad\varepsilon^{(k)}\|^2_{L^2(J;\bL^2(\Omega))},
    &
    \estI &\coloneqq \text{right-hand side of Corollary~\ref{COR:first-apost-Linfty}}, 
    \\
    \errII^2 &\coloneqq \|\grad\varepsilon^{(k)}\|^2_{L^\infty(J;\bL^2(\Omega))}
        + \|\partial_t\varepsilon^{(k)}\|^2_{L^2(J;L^2(\Omega))},
    &
    \estII &\coloneqq \text{right-hand side of Corollary~\ref{COR:second-apost-Linfty}},
\end{aligned}
\end{equation}
together with the effectivity indices $\effI \coloneqq \sqrt{\estI}/\errI$ and $\effII \coloneqq \sqrt{\estII}/\errII$.

For simplicity, we set $\Cdual = \CLtwo = \CHone = 1$ in the elliptic error estimators \eqref{EQ: eta estimator standard}.
The starting values $u_h^0, \dots, u_h^{k-1}$ are the nodal interpolants of the exact solution at $t_0, \dots, t_{k-1}$.
The code is implemented in Julia~\cite{Julia:2017} using the finite element toolbox Ferrite.jl~\cite{Ferrite}.
Since $d=1$, the faces $E \in \Sigma_h$ are the interior nodes, $\norm{\cdot}_{L^2(E)}$ is the absolute value, and $h_E$ denotes the mean of the lengths of the two elements adjacent to $E$.

\subsection{Convergence for $\tau \to 0$}\label{SEC:time-conv}

First, we consider the convergence of the errors and estimators for $\tau \to 0$ using Example~\ref{SEC:ExampleTime}.
We employ uniform \BDF{k} for $k\in \{1, 2, 3, 4\}$ with $\tau=1/N$, $N \in \{64, 92, 128, 180\}$, on a fixed spatial mesh with a mesh size of $h = 1/n_x$, $n_x=16384$, and linear finite elements.

In Table~\ref{TAB:1d-decomp-estI} and Table~\ref{TAB:1d-decomp-estII}, we list the individual error estimator components of Corollary~\ref{COR:first-apost-Linfty} and Corollary~\ref{COR:second-apost-Linfty} without the constants in front of them.
Since the estimators bound the squared error, we expect convergence of order $2k$ for the time-dependent components.
We observe that all temporal error estimators converge with the optimal order.
In total, the error estimator $\NN^{(k)}$ is the dominant term for $k \ge 2$. 
Only for $k=4$ and $N=180$, the spatial error estimator $\Theta^{(4)}(H^1)$ and $\Theta_\infty^{(4)}(H^1)$, which does not decrease with $\tau$, are the largest component.
As expected, the spatial error estimators $\Lambda^{(k)}$, $\Theta_\infty^{(k)}$, and $\Theta^{(k)}$ remain constant since they essentially measure the spatial error.

In Table~\ref{TAB:1d-verify}, we compare the errors $\errI$ and $\errII$ with the estimators $\estI$ and $\estII$ introduced in \eqref{EQ:def-err-est}.
In total, all errors and estimators converge with the optimal order $k$. 
For $k=4$ and $N=180$, saturation due to the fixed spatial mesh occurs. 
The effectivity indices, which are almost independent of $\tau$, grow with $k$.
This is mainly caused by the factor $k$ originating from the discrete Cauchy--Schwarz inequality applied in the proof of Lemma~\ref{LEM:est-residual-final} in front of the estimator $\NN^{(k)}$ in Corollary~\ref{COR:first-apost-Linfty} and Corollary~\ref{COR:second-apost-Linfty}.
Since $\NN^{(k)}$ is the dominant term (see Table~\ref{TAB:1d-decomp-estI} and Table~\ref{TAB:1d-decomp-estII}), the effectivity indices grow with $k$.
As noted in Remark~\ref{REM:factor-k}, this factor $k$ may be optimized.

\begin{table}[h]
\centering
\tiny
\caption{
Example~\ref{SEC:ExampleTime}: Error estimator components of $\mathrm{Est}_I$ in Corollary~\ref{COR:first-apost-Linfty} without the prefactors for \BDF{k} with $\tau=1/N$ and fixed mesh size $h = 1/n_x$, $n_x = 16384$.
}
\label{TAB:1d-decomp-estI}
\resizebox{\textwidth}{!}{%
\begin{tabular}{cccccccccccccc}
\toprule
$k$ & $N$ & \multicolumn{2}{c}{$\FF^{(k)}(H^{-1})$} & \multicolumn{2}{c}{$\NN^{(k)}(H^{-1})$} & \multicolumn{2}{c}{$\Theta^{(k)}(H^1)$} & \multicolumn{2}{c}{$\Theta_\infty^{(k)}(L^2)$} & \multicolumn{2}{c}{$\Lambda^{(k)}(H^{-1})$} & \multicolumn{2}{c}{$\Sest^{(k)}(H^1)$} \\
\midrule
1 & 64 & 5.04e-02 & -- & 0.00e+00 & -- & 1.21e-07 & -- & 4.89e-15 & -- & 2.14e-14 & -- & 6.35e-03 & -- \\
1 & 92 & 2.45e-02 & 1.98 & 0.00e+00 & -- & 1.23e-07 & -0.05 & 5.06e-15 & -0.09 & 2.21e-14 & -0.09 & 3.18e-03 & 1.91 \\
1 & 128 & 1.27e-02 & 1.99 & 0.00e+00 & -- & 1.24e-07 & -0.04 & 5.16e-15 & -0.06 & 2.26e-14 & -0.07 & 1.68e-03 & 1.93 \\
1 & 180 & 6.43e-03 & \textbf{2.00} & 0.00e+00 & \textbf{--} & 1.26e-07 & \textbf{-0.03} & 5.23e-15 & \textbf{-0.04} & 2.30e-14 & \textbf{-0.05} & 8.63e-04 & \textbf{1.95} \\
\addlinespace
2 & 64 & 2.89e-04 & -- & 2.04e-03 & -- & 1.90e-07 & -- & 7.40e-15 & -- & 1.09e-13 & -- & 2.19e-05 & -- \\
2 & 92 & 6.83e-05 & 3.98 & 5.45e-04 & 3.64 & 1.92e-07 & -0.04 & 7.78e-15 & -0.14 & 1.14e-13 & -0.13 & 5.58e-06 & 3.76 \\
2 & 128 & 1.83e-05 & 3.99 & 1.55e-04 & 3.80 & 1.93e-07 & -0.02 & 7.96e-15 & -0.07 & 1.17e-13 & -0.07 & 1.55e-06 & 3.87 \\
2 & 180 & 4.68e-06 & \textbf{3.99} & 4.12e-05 & \textbf{3.89} & 1.94e-07 & \textbf{-0.01} & 8.07e-15 & \textbf{-0.04} & 1.18e-13 & \textbf{-0.04} & 4.07e-07 & \textbf{3.93} \\
\addlinespace
3 & 64 & 2.37e-05 & -- & 4.68e-04 & -- & 2.85e-07 & -- & 1.14e-14 & -- & 3.40e-13 & -- & 1.41e-06 & -- \\
3 & 92 & 2.75e-06 & 5.93 & 6.05e-05 & 5.64 & 2.86e-07 & -0.00 & 1.17e-14 & -0.08 & 3.45e-13 & -0.04 & 1.74e-07 & 5.77 \\
3 & 128 & 3.84e-07 & 5.97 & 8.74e-06 & 5.86 & 2.86e-07 & -0.00 & 1.19e-14 & -0.04 & 3.46e-13 & -0.01 & 2.47e-08 & 5.91 \\
3 & 180 & 4.99e-08 & \textbf{5.98} & 1.15e-06 & \textbf{5.94} & 2.86e-07 & \textbf{0.00} & 1.19e-14 & \textbf{-0.02} & 3.47e-13 & \textbf{-0.00} & 3.23e-09 & \textbf{5.96} \\
\addlinespace
4 & 64 & 3.95e-06 & -- & 1.46e-04 & -- & 4.35e-07 & -- & 1.74e-14 & -- & 8.84e-13 & -- & 2.02e-07 & -- \\
4 & 92 & 2.29e-07 & 7.85 & 8.69e-06 & 7.77 & 4.34e-07 & 0.00 & 1.78e-14 & -0.06 & 8.84e-13 & 0.00 & 1.17e-08 & 7.85 \\
4 & 128 & 1.67e-08 & 7.93 & 6.40e-07 & 7.90 & 4.34e-07 & 0.00 & 1.80e-14 & -0.03 & 8.85e-13 & -0.00 & 8.52e-10 & 7.93 \\
4 & 180 & 1.11e-09 & \textbf{7.96} & 4.27e-08 & \textbf{7.94} & 4.34e-07 & \textbf{0.00} & 1.81e-14 & \textbf{-0.02} & 8.85e-13 & \textbf{-0.00} & 5.65e-11 & \textbf{7.96} \\
\bottomrule
\end{tabular}}
\end{table}

\begin{table}[h]
\centering
\tiny
\caption{
Example~\ref{SEC:ExampleTime}: Error estimator components of $\mathrm{Est}_{II}$ in Corollary~\ref{COR:second-apost-Linfty} without the prefactors for \BDF{k} with $\tau=1/N$ and fixed mesh size $h = 1/n_x$, $n_x = 16384$.}
\label{TAB:1d-decomp-estII}
\resizebox{\textwidth}{!}{%
\begin{tabular}{cccccccccccc}
\toprule
$k$ & $N$ & \multicolumn{2}{c}{$\FF^{(k)}(L^2)$} & \multicolumn{2}{c}{$\NN^{(k)}(L^2)$} & \multicolumn{2}{c}{$\Theta_\infty^{(k)}(H^1)$} & \multicolumn{2}{c}{$\Lambda^{(k)}(L^2)$} & \multicolumn{2}{c}{$\Sest^{(k)}(H^2)$} \\
\midrule
1 & 64 & 4.97e-01 & -- & 0.00e+00 & -- & 1.31e-06 & -- & 2.14e-14 & -- & 6.27e-02 & -- \\
1 & 92 & 2.42e-01 & 1.98 & 0.00e+00 & -- & 1.36e-06 & -0.09 & 2.21e-14 & -0.09 & 3.14e-02 & 1.91 \\
1 & 128 & 1.25e-01 & 1.99 & 0.00e+00 & -- & 1.39e-06 & -0.06 & 2.26e-14 & -0.07 & 1.66e-02 & 1.93 \\
1 & 180 & 6.35e-02 & \textbf{2.00} & 0.00e+00 & \textbf{--} & 1.40e-06 & \textbf{-0.04} & 2.30e-14 & \textbf{-0.05} & 8.52e-03 & \textbf{1.95} \\
\addlinespace
2 & 64 & 2.86e-03 & -- & 2.01e-02 & -- & 1.99e-06 & -- & 1.09e-13 & -- & 2.16e-04 & -- \\
2 & 92 & 6.74e-04 & 3.98 & 5.38e-03 & 3.64 & 2.09e-06 & -0.14 & 1.14e-13 & -0.13 & 5.51e-05 & 3.76 \\
2 & 128 & 1.80e-04 & 3.99 & 1.53e-03 & 3.80 & 2.14e-06 & -0.07 & 1.17e-13 & -0.07 & 1.53e-05 & 3.87 \\
2 & 180 & 4.62e-05 & \textbf{3.99} & 4.06e-04 & \textbf{3.89} & 2.17e-06 & \textbf{-0.04} & 1.18e-13 & \textbf{-0.04} & 4.01e-06 & \textbf{3.93} \\
\addlinespace
3 & 64 & 2.34e-04 & -- & 4.62e-03 & -- & 3.06e-06 & -- & 3.40e-13 & -- & 1.39e-05 & -- \\
3 & 92 & 2.71e-05 & 5.93 & 5.97e-04 & 5.64 & 3.15e-06 & -0.08 & 3.45e-13 & -0.04 & 1.72e-06 & 5.77 \\
3 & 128 & 3.79e-06 & 5.97 & 8.63e-05 & 5.86 & 3.18e-06 & -0.04 & 3.46e-13 & -0.01 & 2.43e-07 & 5.91 \\
3 & 180 & 4.92e-07 & \textbf{5.98} & 1.14e-05 & \textbf{5.94} & 3.20e-06 & \textbf{-0.02} & 3.47e-13 & \textbf{-0.00} & 3.19e-08 & \textbf{5.96} \\
\addlinespace
4 & 64 & 3.90e-05 & -- & 1.44e-03 & -- & 4.67e-06 & -- & 8.84e-13 & -- & 1.99e-06 & -- \\
4 & 92 & 2.26e-06 & 7.85 & 8.58e-05 & 7.77 & 4.78e-06 & -0.06 & 8.84e-13 & 0.00 & 1.15e-07 & 7.85 \\
4 & 128 & 1.65e-07 & 7.93 & 6.32e-06 & 7.90 & 4.83e-06 & -0.03 & 8.85e-13 & -0.00 & 8.41e-09 & 7.93 \\
4 & 180 & 1.09e-08 & \textbf{7.96} & 4.21e-07 & \textbf{7.94} & 4.86e-06 & \textbf{-0.02} & 8.85e-13 & \textbf{-0.00} & 5.58e-10 & \textbf{7.96} \\
\bottomrule
\end{tabular}}
\end{table}

\begin{table}[htb!]
\centering
\footnotesize
\caption{Example~\ref{SEC:ExampleTime}: errors $\mathcal E_I$, $\mathcal E_{II}$ and estimators $\sqrt{\mathrm{Est}_I}$, $\sqrt{\mathrm{Est}_{II}}$ for \BDF{k} with $\tau=1/N$ on a fixed mesh with $n_x=16384$ and $\mathbb P_1$ elements.}
\label{TAB:1d-verify}
\setlength{\tabcolsep}{4pt}
\begin{tabular}{cc ccccc ccccc}
\toprule
 & & \multicolumn{5}{c}{Corollary~\ref{COR:first-apost-Linfty}} & \multicolumn{5}{c}{Corollary~\ref{COR:second-apost-Linfty}} \\
\cmidrule(lr){3-7}\cmidrule(lr){8-12}
$k$ & $N$ & \multicolumn{2}{c}{$\mathcal E_I$} & \multicolumn{2}{c}{$\sqrt{\mathrm{Est}_I}$} & $\mathrm{Eff}_I$ & \multicolumn{2}{c}{$\mathcal E_{II}$} & \multicolumn{2}{c}{$\sqrt{\mathrm{Est}_{II}}$} & $\mathrm{Eff}_{II}$ \\
\midrule
1 & 64 & 9.25e-02 & -- & 8.25e-01 & -- & 8.92 & 6.32e-01 & -- & 2.59e+00 & -- & 4.10 \\
1 & 92 & 6.54e-02 & 0.95 & 5.77e-01 & 0.99 & 8.82 & 4.47e-01 & 0.95 & 1.81e+00 & 0.99 & 4.05 \\
1 & 128 & 4.75e-02 & 0.97 & 4.15e-01 & 0.99 & 8.75 & 3.25e-01 & 0.97 & 1.31e+00 & 0.99 & 4.02 \\
1 & 180 & 3.39e-02 & \textbf{0.99} & 2.96e-01 & \textbf{1.00} & 8.72 & 2.33e-01 & \textbf{0.98} & 9.30e-01 & \textbf{1.00} & 4.00 \\
\addlinespace
2 & 64 & 2.48e-02 & -- & 2.37e-01 & -- & 9.57 & 1.86e-01 & -- & 7.45e-01 & -- & 4.00 \\
2 & 92 & 1.26e-02 & 1.87 & 1.22e-01 & 1.84 & 9.68 & 9.55e-02 & 1.84 & 3.82e-01 & 1.84 & 4.00 \\
2 & 128 & 6.63e-03 & 1.93 & 6.47e-02 & 1.91 & 9.75 & 5.07e-02 & 1.92 & 2.03e-01 & 1.91 & 4.01 \\
2 & 180 & 3.38e-03 & \textbf{1.97} & 3.33e-02 & \textbf{1.95} & 9.83 & 2.60e-02 & \textbf{1.96} & 1.04e-01 & \textbf{1.95} & 4.01 \\
\addlinespace
3 & 64 & 9.45e-03 & -- & 1.33e-01 & -- & 14.09 & 8.17e-02 & -- & 4.18e-01 & -- & 5.12 \\
3 & 92 & 3.36e-03 & 2.85 & 4.77e-02 & 2.83 & 14.21 & 2.92e-02 & 2.84 & 1.50e-01 & 2.83 & 5.14 \\
3 & 128 & 1.26e-03 & 2.96 & 1.81e-02 & 2.93 & 14.36 & 1.10e-02 & 2.94 & 5.70e-02 & 2.93 & 5.16 \\
3 & 180 & 4.59e-04 & \textbf{2.97} & 6.63e-03 & \textbf{2.95} & 14.45 & 4.01e-03 & \textbf{2.98} & 2.08e-02 & \textbf{2.95} & 5.20 \\
\addlinespace
4 & 64 & 4.63e-03 & -- & 8.47e-02 & -- & 18.29 & 4.38e-02 & -- & 2.66e-01 & -- & 6.07 \\
4 & 92 & 1.13e-03 & 3.89 & 2.07e-02 & 3.88 & 18.32 & 1.07e-02 & 3.89 & 6.51e-02 & 3.88 & 6.09 \\
4 & 128 & 3.09e-04 & 3.93 & 5.69e-03 & 3.91 & 18.42 & 2.90e-03 & 3.95 & 1.79e-02 & 3.91 & 6.19 \\
4 & 180 & 8.71e-05 & \textbf{3.72} & 1.72e-03 & \textbf{3.50} & 19.80 & 7.56e-04 & \textbf{3.94} & 5.52e-03 & \textbf{3.45} & 7.30 \\
\bottomrule
\end{tabular}
\end{table}

\subsection{Convergence for $h \to 0$}

Second, we consider the convergence of the errors and estimators for $h \to 0$ using Example~\ref{SEC:ExampleSpace}.
We employ uniform \BDF{2} with $\tau=1/32$ and various spatial meshes with finite element degree $p \in \{1, 2\}$.

In Table~\ref{TAB:1d-smooth-h-decomp-estI} and Table~\ref{TAB:1d-smooth-h-decomp-estII}, we present the individual error estimator components. 
Again, since the estimators bound the squared error, we expect convergence of (at least) order $2p$ for the spatial error estimators, which we observe for each estimator. 
Together with Table~\ref{TAB:1d-decomp-estI} and Table~\ref{TAB:1d-decomp-estII}, we observe that we effectively split the error estimators into time and space contributions, and each estimator converges either in time or space with the optimal order.

In Table~\ref{TAB:1d-smooth-h-est}, we compare the errors $\errI$ and $\errII$ with the total estimators $\estI$ and $\estII$ introduced in \eqref{EQ:def-err-est}.
All errors and estimators converge with the optimal order $p$.
The effectivity indices are larger than those for the convergence in time.
They result from the overstimation by the residual-type estimator $\eta$, whose constants we set for simplicity to $1$,~cf.~Section~\ref{SEC:impl-details}, ans from the constants in $\Theta^{(k)}$ and $\Theta_\infty^{(k)}$.
Such an overestimation is known for residual-type error estimators,~cf.~\cite[Fig.~1]{CarstensenBartelsKlose:2001} and \cite[Fig.~1.1]{CarstensenMerdon:2010}.
Since Assumption~\ref{ASS:a-post-bound-ell-rec} only requires computable bounds for the elliptic reconstruction, the estimator $\eta$ can be replaced by any other reliable estimator for the elliptic problem.

\begin{table}[h]
\centering
\tiny
\caption{
Example~\ref{SEC:ExampleSpace}: Error estimator components of $\mathrm{Est}_I$ in Corollary~\ref{COR:first-apost-Linfty} without the prefactors for \BDF{2} with fixed $\tau=1/32$ and various mesh size $h = 1/n_x$ with finite element polynomial degree $p$.
}
\label{TAB:1d-smooth-h-decomp-estI}
\resizebox{\textwidth}{!}{%
\begin{tabular}{cccccccccccccc}
\toprule
$\mathbb P$ & $n_x$ & \multicolumn{2}{c}{$\FF^{(2)}(H^{-1})$} & \multicolumn{2}{c}{$\NN^{(2)}(H^{-1})$} & \multicolumn{2}{c}{$\Theta^{(2)}(H^1)$} & \multicolumn{2}{c}{$\Theta_\infty^{(2)}(L^2)$} & \multicolumn{2}{c}{$\Lambda^{(2)}(H^{-1})$} & \multicolumn{2}{c}{$\Sest^{(2)}(H^1)$} \\
\midrule
1 & 8 & 1.73e-08 & -- & 4.49e-08 & -- & 3.65e+00 & -- & 1.29e-01 & -- & 4.76e-02 & -- & 2.11e-08 & -- \\
1 & 16 & 1.30e-08 & 0.42 & 3.70e-09 & 3.60 & 9.23e-01 & 1.98 & 8.17e-03 & 3.99 & 3.00e-03 & 3.99 & 1.75e-08 & 0.27 \\
1 & 32 & 1.27e-08 & 0.03 & 1.55e-09 & 1.26 & 2.31e-01 & 2.00 & 5.12e-04 & 4.00 & 1.88e-04 & 4.00 & 1.66e-08 & 0.08 \\
1 & 64 & 1.27e-08 & \textbf{0.00} & 1.55e-09 & \textbf{0.00} & 5.79e-02 & \textbf{2.00} & 3.20e-05 & \textbf{4.00} & 1.17e-05 & \textbf{4.00} & 1.63e-08 & \textbf{0.02} \\
\addlinespace
2 & 8 & 1.27e-08 & -- & 1.81e-09 & -- & 2.38e-02 & -- & 8.39e-04 & -- & 4.84e-06 & -- & 1.72e-08 & -- \\
2 & 16 & 1.27e-08 & -0.00 & 1.59e-09 & 0.18 & 1.49e-03 & 4.00 & 1.31e-05 & 6.00 & 1.89e-08 & 8.00 & 1.63e-08 & 0.08 \\
2 & 32 & 1.27e-08 & -0.00 & 1.59e-09 & 0.00 & 9.30e-05 & 4.00 & 2.06e-07 & 6.00 & 7.37e-11 & 8.00 & 1.63e-08 & 0.00 \\
2 & 64 & 1.27e-08 & \textbf{-0.00} & 1.59e-09 & \textbf{-0.00} & 5.81e-06 & \textbf{4.00} & 3.21e-09 & \textbf{6.00} & 2.88e-13 & \textbf{8.00} & 1.63e-08 & \textbf{0.00} \\
\bottomrule
\end{tabular}}
\end{table}

\begin{table}[h]
\centering
\tiny
\caption{
Example~\ref{SEC:ExampleSpace}: Error estimator components of $\mathrm{Est}_{II}$ in Corollary~\ref{COR:second-apost-Linfty} without the prefactors for \BDF{2} with fixed $\tau=1/32$ and various mesh size $h = 1/n_x$ with finite element polynomial degree $p$.
}
\label{TAB:1d-smooth-h-decomp-estII}
\resizebox{\textwidth}{!}{%
\begin{tabular}{cccccccccccc}
\toprule
$\mathbb P$ & $n_x$ & \multicolumn{2}{c}{$\FF^{(2)}(L^2)$} & \multicolumn{2}{c}{$\NN^{(2)}(L^2)$} & \multicolumn{2}{c}{$\Theta_\infty^{(2)}(H^1)$} & \multicolumn{2}{c}{$\Lambda^{(2)}(L^2)$} & \multicolumn{2}{c}{$\Sest^{(2)}(H^2)$} \\
\midrule
1 & 8 & 5.42e-04 & -- & 4.49e-07 & -- & 8.28e+00 & -- & 4.76e-02 & -- & 1.61e-07 & -- \\
1 & 16 & 3.31e-05 & 4.03 & 3.66e-08 & 3.62 & 2.09e+00 & 1.99 & 3.00e-03 & 3.99 & 1.61e-07 & 0.00 \\
1 & 32 & 2.17e-06 & 3.93 & 1.53e-08 & 1.26 & 5.24e-01 & 2.00 & 1.88e-04 & 4.00 & 1.61e-07 & -0.00 \\
1 & 64 & 2.53e-07 & \textbf{3.10} & 1.53e-08 & \textbf{0.00} & 1.31e-01 & \textbf{2.00} & 1.17e-05 & \textbf{4.00} & 1.61e-07 & \textbf{-0.00} \\
\addlinespace
2 & 8 & 1.95e-06 & -- & 9.19e-07 & -- & 5.37e-02 & -- & 3.09e-04 & -- & 2.03e-07 & -- \\
2 & 16 & 1.55e-07 & 3.65 & 7.59e-08 & 3.60 & 3.37e-03 & 4.00 & 4.83e-06 & 6.00 & 1.63e-07 & 0.31 \\
2 & 32 & 1.26e-07 & 0.30 & 1.95e-08 & 1.96 & 2.11e-04 & 4.00 & 7.55e-08 & 6.00 & 1.61e-07 & 0.02 \\
2 & 64 & 1.26e-07 & \textbf{0.01} & 1.60e-08 & \textbf{0.29} & 1.32e-05 & \textbf{4.00} & 1.18e-09 & \textbf{6.00} & 1.61e-07 & \textbf{0.00} \\
\bottomrule
\end{tabular}}
\end{table}

\begin{table}[htb!]
\centering
\footnotesize
\caption{Example~\ref{SEC:ExampleSpace}: errors $\mathcal E_I$, $\mathcal E_{II}$ and estimators $\sqrt{\mathrm{Est}_I}$, $\sqrt{\mathrm{Est}_{II}}$ of Corollaries~\ref{COR:first-apost-Linfty} and~\ref{COR:second-apost-Linfty} for \BDF{2} with $\tau=1/32$ and $\mathbb P_p$ elements on uniform meshes with $h=1/n_x$.}
\label{TAB:1d-smooth-h-est}
\setlength{\tabcolsep}{4pt}
\begin{tabular}{cc ccccc ccccc}
\toprule
 & & \multicolumn{5}{c}{Corollary~\ref{COR:first-apost-Linfty}} & \multicolumn{5}{c}{Corollary~\ref{COR:second-apost-Linfty}} \\
\cmidrule(lr){3-7}\cmidrule(lr){8-12}
$\mathbb P$ & $n_x$ & \multicolumn{2}{c}{$\mathcal E_I$} & \multicolumn{2}{c}{$\sqrt{\mathrm{Est}_I}$} & $\mathrm{Eff}_I$ & \multicolumn{2}{c}{$\mathcal E_{II}$} & \multicolumn{2}{c}{$\sqrt{\mathrm{Est}_{II}}$} & $\mathrm{Eff}_{II}$ \\
\midrule
1 & 8 & 1.60e-01 & -- & 2.88e+00 & -- & 18.06 & 2.44e-01 & -- & 5.39e+00 & -- & 22.12 \\
1 & 16 & 7.98e-02 & 1.00 & 1.38e+00 & 1.06 & 17.31 & 1.22e-01 & 1.00 & 2.68e+00 & 1.01 & 21.99 \\
1 & 32 & 3.99e-02 & 1.00 & 6.83e-01 & 1.02 & 17.11 & 6.10e-02 & 1.00 & 1.34e+00 & 1.00 & 21.96 \\
1 & 64 & 2.00e-02 & \textbf{1.00} & 3.41e-01 & \textbf{1.00} & 17.06 & 3.05e-02 & \textbf{1.00} & 6.70e-01 & \textbf{1.00} & 21.95 \\
\addlinespace
2 & 8 & 8.08e-03 & -- & 2.25e-01 & -- & 27.78 & 1.23e-02 & -- & 4.31e-01 & -- & 34.92 \\
2 & 16 & 2.02e-03 & 2.00 & 5.49e-02 & 2.03 & 27.15 & 3.09e-03 & 2.00 & 1.07e-01 & 2.01 & 34.62 \\
2 & 32 & 5.08e-04 & 1.99 & 1.37e-02 & 2.01 & 26.92 & 7.77e-04 & 1.99 & 2.68e-02 & 2.00 & 34.50 \\
2 & 64 & 1.35e-04 & \textbf{1.91} & 3.49e-03 & \textbf{1.97} & 25.76 & 2.09e-04 & \textbf{1.90} & 7.06e-03 & \textbf{1.92} & 33.84 \\
\bottomrule
\end{tabular}
\end{table}

\subsection{Variable time-steps} 

Finally, we consider \VSBDF{k}, $k \in \{1, 2, 3\}$, with variable time-steps for Example~\ref{SEC:ExampleTime}. 
For $N \in \{64, 92, 128, 180\}$, we use $N$ time-steps alternating between $\tau$ and $2 \tau$ with $\tau = \frac{2T}{3N}$.
Thus, the step-ratios $\kappa_n = \tau_n / \tau_{n-1}$ alternate between $2$ and $\tfrac{1}{2}$.
The spatial discretization is the same as in Section~\ref{SEC:time-conv}.
Table~\ref{TAB:1d-verify-variable} shows that the errors and the estimators converge with the optimal order $k$ also for variable time-steps. 
The effectivity indices are (almost) independent of $N$ and of the same size as for uniform time-steps,~cf.~Table~\ref{TAB:1d-verify}.
For $k=2, 3$, the effectivity indices are even slightly smaller.

\begin{table}[htb!]
\centering
\footnotesize
\caption{Example~\ref{SEC:ExampleTime}: errors $\mathcal E_I$, $\mathcal E_{II}$ and estimators $\sqrt{\mathrm{Est}_I}$, $\sqrt{\mathrm{Est}_{II}}$ for \VSBDF{k} using alternating time-steps $\tau$ and $2\tau$ with $\tau = \frac{2T}{3N}$ on a fixed mesh with $n_x=16384$ and $\mathbb P_1$ elements.}
\label{TAB:1d-verify-variable}
\setlength{\tabcolsep}{4pt}
\begin{tabular}{cc ccccc ccccc}
\toprule
 & & \multicolumn{5}{c}{Corollary~\ref{COR:first-apost-Linfty}} & \multicolumn{5}{c}{Corollary~\ref{COR:second-apost-Linfty}} \\
\cmidrule(lr){3-7}\cmidrule(lr){8-12}
$k$ & $N$ & \multicolumn{2}{c}{$\mathcal E_I$} & \multicolumn{2}{c}{$\sqrt{\mathrm{Est}_I}$} & $\mathrm{Eff}_I$ & \multicolumn{2}{c}{$\mathcal E_{II}$} & \multicolumn{2}{c}{$\sqrt{\mathrm{Est}_{II}}$} & $\mathrm{Eff}_{II}$ \\
\midrule
1 & 64 & 1.03e-01 & -- & 9.49e-01 & -- & 9.19 & 7.20e-01 & -- & 2.98e+00 & -- & 4.14 \\
1 & 92 & 7.23e-02 & 0.98 & 6.64e-01 & 0.98 & 9.19 & 5.11e-01 & 0.95 & 2.09e+00 & 0.98 & 4.09 \\
1 & 128 & 5.24e-02 & 0.98 & 4.79e-01 & 0.99 & 9.15 & 3.72e-01 & 0.96 & 1.51e+00 & 0.99 & 4.05 \\
1 & 180 & 3.76e-02 & \textbf{0.98} & 3.41e-01 & \textbf{0.99} & 9.09 & 2.67e-01 & \textbf{0.97} & 1.07e+00 & \textbf{0.99} & 4.02 \\
\addlinespace
2 & 64 & 2.52e-02 & -- & 2.33e-01 & -- & 9.28 & 1.88e-01 & -- & 7.33e-01 & -- & 3.90 \\
2 & 92 & 1.27e-02 & 1.89 & 1.19e-01 & 1.86 & 9.37 & 9.60e-02 & 1.85 & 3.73e-01 & 1.86 & 3.88 \\
2 & 128 & 6.65e-03 & 1.95 & 6.28e-02 & 1.93 & 9.45 & 5.09e-02 & 1.92 & 1.97e-01 & 1.93 & 3.88 \\
2 & 180 & 3.39e-03 & \textbf{1.98} & 3.22e-02 & \textbf{1.96} & 9.50 & 2.61e-02 & \textbf{1.96} & 1.01e-01 & \textbf{1.96} & 3.88 \\
\addlinespace
3 & 64 & 9.40e-03 & -- & 1.17e-01 & -- & 12.48 & 8.20e-02 & -- & 3.69e-01 & -- & 4.50 \\
3 & 92 & 3.37e-03 & 2.83 & 4.15e-02 & 2.86 & 12.34 & 2.92e-02 & 2.85 & 1.30e-01 & 2.86 & 4.47 \\
3 & 128 & 1.26e-03 & 2.99 & 1.57e-02 & 2.94 & 12.51 & 1.10e-02 & 2.95 & 4.93e-02 & 2.94 & 4.48 \\
3 & 180 & 4.57e-04 & \textbf{2.96} & 5.74e-03 & \textbf{2.95} & 12.56 & 3.99e-03 & \textbf{2.98} & 1.81e-02 & \textbf{2.94} & 4.54 \\
\bottomrule
\end{tabular}
\end{table}

\appendix
\section{Computing the remainder polynomials}\label{AP:remaider-poly}

In this appendix, we provide the proof of Lemma~\ref{LEM:rho_i}.

\begin{proof}[Proof of Lemma~\ref{LEM:rho_i}]
    To obtain the closed-form \eqref{EQ:rho-closed-form}, it is convenient to derive $\rho_k^n$ recursively. 
    We start with the known recurrence relation for the Newton polynomials
    \begin{align*}
        \omega_k^{n}(t) = \omega_{k-1}^{n}(t)(t-t_{n-k+1}).
    \end{align*}
    Differentiating this product yields
    \begin{align*}
        \ptl_t \omega_k^{n}(t) = \ptl_t \omega_{k-1}^{n}(t)(t-t_{n-k+1}) + \omega_{k-1}^{n}(t).
    \end{align*}
    Evaluating at $t_n$ and noting that $\omega_{k-1}^{n}(t_n) = 0$ for $k \ge 2$, we find
    \begin{align*}
        \ptl_t \omega_k^{n}(t_n) = \ptl_t \omega_{k-1}^{n}(t_n)(t_n-t_{n-k+1}).
    \end{align*}
    Using the definition \eqref{EQ:Def_rho},~i.e.,~$\rho_k^n(t) = \ptl_t \omega_k^{n}(t) - \ptl_t \omega_k^{n}(t_n)$, substituting the recursive formulas and grouping the terms associated with $\rho_{k-1}^n(t)$ we obtain the recursive relation
    \begin{align*}
        \rho_k^n(t) = \rho_{k-1}^n(t)(t-t_{n-k+1}) + \ptl_t \omega_{k-1}^{n}(t_n)(t-t_n) + \omega_{k-1}^{n}(t),
    \end{align*}
    which holds for $k \ge 2$ with the initial condition $\rho_1^n(t) = 0$.
    
    \textbf{\boldmath $k=2$:}    
    \begin{align*}
        \rho_2^n(t) &= \rho_1^n(t)(t-t_{n-1}) + \ptl_t \omega_1^{n}(t_n)(t-t_n) + \omega_1^{n}(t) \\
                  &= 0 \cdot (t-t_{n-1}) + 1 \cdot (t-t_n) + (t-t_n) \\
                  &= 2(t-t_n).
    \end{align*}

    \textbf{\boldmath $k=3$:}    
    We use $\rho_2^n(t) = 2(t-t_n)$ and $\omega_2^{n}(t) = (t-t_n)((t-t_n) + c_1^n)$. The required derivative evaluation is $\ptl_t \omega_2^{n}(t_n) = c_1^n$. We also rewrite the shifted root as $(t-t_{n-2}) = (t-t_n) + c_2^n$. Inserting these into the recursion yields
    \begin{align*}
        \rho_3^n(t) &= \rho_2^n(t)(t-t_{n-2}) + \ptl_t \omega_2^{n}(t_n)(t-t_n) + \omega_2^{n}(t) \\
                  &= 2(t-t_n) \big( (t-t_n) + c_2^n \big) + c_1^{n}(t-t_n) + (t-t_n) \big( (t-t_n) + c_1^n \big).
    \end{align*}
    We expand the products and group the coefficients to obtain
    \begin{align*}
        \rho_3^n(t) &= 3(t-t_n)^2 + 2(c_1^n + c_2^n)(t-t_n).
    \end{align*}

    \textbf{\boldmath $k=4$:}    
    We have $\omega_3^{n}(t) = (t-t_n)^3 + (c_1^n+c_2^n)(t-t_n)^2 + c_1^n c_2^n (t-t_n)$ and the evaluation of its derivative provides $\ptl_t \omega_3^{n}(t_n) = c_1^n c_2^n$. 
    We express the final root as $(t-t_{n-3}) = (t-t_n) + c_3^n$. 
    Expanding
    \begin{align*}
        \rho_3^n(t)(t-t_{n-3}) &= 3(t-t_n)^3 + 3c_3^n(t-t_n)^2 + 2(c_1^n+c_2^n)(t-t_n)^2 + 2c_3^n(c_1^n+c_2^n)(t-t_n).
    \end{align*}
    and applying the recursive relation gives
    \begin{align*}
        \rho_4^n(t) &= \rho_3^n(t)(t-t_{n-3}) + \ptl_t \omega_3^{n}(t_n)(t-t_n) + \omega_3^{n}(t) \\
                  &= \big[ 3(t-t_n)^2 + 2(c_1^n + c_2^n)(t-t_n) \big] \big( (t-t_n) + c_3^n \big) + c_1^n c_2^n (t-t_n) \\
                  &\quad + \big[ (t-t_n)^3 + (c_1^n+c_2^n)(t-t_n)^2 + c_1^n c_2^n (t-t_n) \big].
    \end{align*}   
    Next, we sum the coefficients for each polynomial degree. 
    For the cubic term, we obtain $3(t-t_n)^3 + (t-t_n)^3 = 4(t-t_n)^3$.
    For the quadratic term, we group $3c_3^n + 2(c_1^n+c_2^n) + (c_1^n+c_2^n) = 3(c_1^n+c_2^n+c_3^n)$.
    For the linear term, we add $2(c_1^n c_3^n + c_2^n c_3^n) + c_1^n c_2^n + c_1^n c_2^n = 2(c_1^n c_2^n + c_1^n c_3^n + c_2^n c_3^n)$.
    Combining these degree-wise summations yields 
    \begin{align*}
        \rho_4^n(t) &= 4(t-t_n)^3 + 3(c_1^n+c_2^n+c_3^n)(t-t_n)^2 + 2(c_1^n c_2^n + c_1^n c_3^n + c_2^n c_3^n)(t-t_n).
    \end{align*}
\end{proof}

\section{Approximation at negative times}\label{APP:approx-at-ge-times}
Using \eqref{EQ:vhn-init-def}, \eqref{EQ:def-neg-time-vals-u} and \eqref{EQ:def-delta-coeff}, we have 
\begin{align*}
    v_h^{n, (k)} = \Delta_h u_h^n + f_h^n = \sum\limits_{i=1}^k \ptl_t \omega_i^{n}(t_n) \nabla^i u_h^n
    = 
    \sum\limits_{i=1}^k \ptl_t \omega_i^{n}(t_n) \sum\limits_{j=0}^i \delta_j^{i, n} u_h^{n-j}
\end{align*}
for $1 \le n \le k-1$.
Defining $\mu_j^{n,(k)} \coloneqq \sum_{i=j}^{k} \ptl_t \omega_i^{n}(t_n)\,\delta_j^{i,n}$ and changing the order of summation gives
\begin{align*}
    v_h^{n, (k)} 
    = 
    \sum\limits_{j=0}^k \mu_j^{n,(k)} u_h^{n-j}
\end{align*}
for $1 \le n \le k-1$.
Thus, we obtain the values at negative times $u_h^{-\ell}$ recursively for $\ell=1$ to $k-1$ by 
\begin{align*}
    u_h^{-\ell} = p\big(v_h^{k - \ell, (k)}, u_h^{k - \ell},\dots,u_h^{1 - \ell}\big)
    \coloneqq \frac{1}{\mu_k^{k - \ell,(k)}}\Big(v_h^{k - \ell, (k)} - \sum_{j=0}^{k-1}\mu_j^{k - \ell,(k)} u_h^{k - \ell -j}\Big).
\end{align*}
This is well-defined, since $\delta_k^{k, n} = (-1)^k / \prod_{m=1}^k c_m^{n-k+m}$ by the recursion \eqref{EQ:def-delta-coeff-recursion}, and thus
\begin{align*}
    \mu_k^{n, (k)} = \ptl_t \omega_k^n(t_n) \delta_k^{k, n} = (-1)^k \frac{\prod_{m=1}^{k-1} c_m^n}{ \prod_{m=1}^k c_m^{n-k+m}} \neq 0.
\end{align*}

\section{Integrating Newton polynomials}
The following integrals are used in the proof of Lemma~\ref{LEM:est-residual-final}.

\begin{Lemma}\label{LEM:omega-moments}
    For $\omega_i^n$ defined in \eqref{EQ:Newton-interpol-poly}, it holds
    \begin{align}
        \label{EQ:mom-11}
        \int_{I_n}\omega_1^n(t)^2\,\dt &= \frac{\tau_n^3}{3},
        \\
        \label{EQ:mom-12}
        \int_{I_n}\omega_1^n(t)\,\omega_2^n(t)\,\dt &= \frac{\tau_n^4}{12},
        \\
        \label{EQ:mom-22}
        \int_{I_n}\omega_2^n(t)^2\,\dt &= \frac{\tau_n^5}{30}.
    \end{align}
\end{Lemma}
\begin{proof}
    First, we substitute $s=t-t_n\in(-\tau_n,0)$ for $t \in I_n$. 
    Using $\omega_1^n(t)=s$ and $\omega_2^n(t)=s(s+c_1^n)$ with $c_1^n=\tau_n$, direct integration yields
    \begin{align*}
        \int_{-\tau_n}^0 s^2 \dss &= \frac{\tau_n^3}{3},
        \qq 
        \int_{-\tau_n}^0 s^2(s+\tau_n) \dss 
        = \frac{\tau_n^4}{12},
        \qq
        \int_{-\tau_n}^0 s^2(s+\tau_n)^2 \dss
        = \frac{\tau_n^5}{30}.
    \end{align*}
\end{proof}

{\sloppy
\subsection*{Acknowledgement}
The author gratefully acknowledges the support of the Deutsche Forschungsgemeinschaft (DFG, German Research Foundation) within the SFB 1173 \enquote{Wave Phenomena} (Project ID 258734477), as well as the ERC grant OPTIMAL (grant agreement No. 101125225).
\par}

\subsection*{Declaration}
The author used Claude Opus 5 (Anthropic) for proofreading, and reviewing the code for the numerical experiments.
The author takes full responsibility for the content of this work.

\bibliographystyle{abbrv}
\bibliography{references}

\begin{thebibliography}{10}

\bibitem{AkrivisChatzipantelidis:2010}
G.~Akrivis and P.~Chatzipantelidis.
\newblock A posteriori error estimates for the two-step backward
  differentiation formula method for parabolic equations.
\newblock {\em SIAM J. Numer. Anal.}, 48(1):109--132, 2010.

\bibitem{AkrivisET:2021}
G.~Akrivis, M.~Chen, F.~Yu, and Z.~Zhou.
\newblock The energy technique for the six-step {BDF} method.
\newblock {\em SIAM J. Numer. Anal.}, 59(5):2449--2472, 2021.

\bibitem{AkrivisLubich:2015}
G.~Akrivis and C.~Lubich.
\newblock Fully implicit, linearly implicit and implicit-explicit backward
  difference formulae for quasi-linear parabolic equations.
\newblock {\em Numer. Math.}, 131(4):713--735, 2015.

\bibitem{AkrivisMakridakisNochetto:2006}
G.~Akrivis, C.~G. Makridakis, and R.~H. Nochetto.
\newblock A posteriori error estimates for the {C}rank--{N}icolson method for
  parabolic equations.
\newblock {\em Math. Comp.}, 75(254):511--531, 2006.

\bibitem{AkrivisMakridakisNochetto:2009}
G.~Akrivis, C.~G. Makridakis, and R.~H. Nochetto.
\newblock Optimal order a posteriori error estimates for a class of
  {R}unge--{K}utta and {G}alerkin methods.
\newblock {\em Numer. Math.}, 114(1):133--160, 2009.

\bibitem{AkrivisMakridakisNochetto:2011}
G.~Akrivis, C.~G. Makridakis, and R.~H. Nochetto.
\newblock Galerkin and {R}unge--{K}utta methods: unified formulation, a
  posteriori error estimates and nodal superconvergence.
\newblock {\em Numer. Math.}, 118(3):429--456, 2011.

\bibitem{BaenschBrenner:2019}
E.~B\"ansch and A.~Brenner.
\newblock A posteriori estimates for the two-step backward differentiation
  formula and discrete regularity for the time-dependent {S}tokes equations.
\newblock {\em IMA J. Numer. Anal.}, 39(2):713--759, 2019.

\bibitem{BaenschKarakatsaniMakridakis:2012}
E.~B\"ansch, F.~Karakatsani, and C.~G. Makridakis.
\newblock A posteriori error control for fully discrete {C}rank--{N}icolson
  schemes.
\newblock {\em SIAM J. Numer. Anal.}, 50(6):2845--2872, 2012.

\bibitem{Becker:1998}
J.~Becker.
\newblock A second order backward difference method with variable steps for a
  parabolic problem.
\newblock {\em BIT}, 38(4):644--662, 1998.

\bibitem{Julia:2017}
J.~Bezanson, A.~Edelman, S.~Karpinski, and V.~B. Shah.
\newblock Julia: A fresh approach to numerical computing.
\newblock {\em SIAM Rev.}, 59(1):65--98, 2017.

\bibitem{CalvoGrandeGrigorieff:1990}
M.~Calvo, T.~Grande, and R.~D. Grigorieff.
\newblock On the zero stability of the variable order variable stepsize
  {BDF}-formulas.
\newblock {\em Numer. Math.}, 57(1):39--50, 1990.

\bibitem{Ferrite}
K.~Carlsson, F.~Ekre, and {Ferrite.jl contributors}.
\newblock {Ferrite.jl}.
\newblock \url{https://github.com/Ferrite-FEM/Ferrite.jl}, 2026.
\newblock Version 1.5.0, \doi{10.5281/zenodo.13862652}.

\bibitem{CarstensenBartelsKlose:2001}
C.~Carstensen, S.~Bartels, and R.~Klose.
\newblock An experimental survey of a posteriori {C}ourant finite element error
  control for the {P}oisson equation.
\newblock {\em Adv. Comput. Math.}, 15(1--4):79--106, 2001.

\bibitem{CarstensenMerdon:2010}
C.~Carstensen and C.~Merdon.
\newblock Estimator competition for {P}oisson problems.
\newblock {\em J. Comput. Math.}, 28(3):309--330, 2010.

\bibitem{CrouzeixLisbona:1984}
M.~Crouzeix and F.~J. Lisbona.
\newblock The convergence of variable-stepsize, variable-formula, multistep
  methods.
\newblock {\em SIAM J. Numer. Anal.}, 21(3):512--534, 1984.

\bibitem{Emmrich:2005}
E.~Emmrich.
\newblock Stability and error of the variable two-step {BDF} for semilinear
  parabolic problems.
\newblock {\em J. Appl. Math. Comput.}, 19(1--2):33--55, 2005.

\bibitem{ErikssonJohnson:1991}
K.~Eriksson and C.~Johnson.
\newblock Adaptive finite element methods for parabolic problems. {I}. {A}
  linear model problem.
\newblock {\em SIAM J. Numer. Anal.}, 28(1):43--77, 1991.

\bibitem{ErnSmearsVohralik:2017}
A.~Ern, I.~Smears, and M.~Vohral\'{\i}k.
\newblock Guaranteed, locally space-time efficient, and polynomial-degree
  robust a~posteriori error estimates for high-order discretizations of
  parabolic problems.
\newblock {\em SIAM J. Numer. Anal.}, 55(6):2811--2834, 2017.

\bibitem{ErnVohralik:2010}
A.~Ern and M.~Vohral\'{\i}k.
\newblock A posteriori error estimation based on potential and flux
  reconstruction for the heat equation.
\newblock {\em SIAM J. Numer. Anal.}, 48(1):198--223, 2010.

\bibitem{Evans:1998}
L.~C. Evans.
\newblock {\em Partial differential equations}, volume~19 of {\em Graduate
  Studies in Mathematics}.
\newblock American Mathematical Society, Providence, RI, 1998.

\bibitem{GeorgoulisMakridakis:2023}
E.~H. Georgoulis and C.~G. Makridakis.
\newblock Lower bounds, elliptic reconstruction and a posteriori error control
  of parabolic problems.
\newblock {\em IMA J. Numer. Anal.}, 43(6):3212--3242, 2023.

\bibitem{Grigorieff:1983}
R.~D. Grigorieff.
\newblock Stability of multistep-methods on variable grids.
\newblock {\em Numer. Math.}, 42(3):359--377, 1983.

\bibitem{HairerWanner:1996}
E.~Hairer and G.~Wanner.
\newblock {\em Solving ordinary differential equations {II}: Stiff and
  differential-algebraic problems}, volume~14 of {\em Springer Series in
  Computational Mathematics}.
\newblock Springer-Verlag, Berlin, second edition, 1996.

\bibitem{Karch:2026}
S.~Karch.
\newblock A posteriori error analysis for the second-order {BDF} method for the
  {L}andau--{L}ifshitz--{G}ilbert equation.
\newblock arXiv:2608.02021, 2026.

\bibitem{LakkisMakridakis:2006}
O.~Lakkis and C.~G. Makridakis.
\newblock Elliptic reconstruction and a posteriori error estimates for fully
  discrete linear parabolic problems.
\newblock {\em Math. Comp.}, 75(256):1627--1658, 2006.

\bibitem{LiLiao:2022}
Z.~Li and H.-L. Liao.
\newblock Stability of variable-step {BDF}2 and {BDF}3 methods.
\newblock {\em SIAM J. Numer. Anal.}, 60(4):2253--2272, 2022.

\bibitem{LiaoTangZhou:2023}
H.-L. Liao, T.~Tang, and T.~Zhou.
\newblock Discrete energy analysis of the third-order variable-step {BDF}
  time-stepping for diffusion equations.
\newblock {\em J. Comput. Math.}, 41(2):325--344, 2023.

\bibitem{LiaoZhang:2021}
H.-L. Liao and Z.~Zhang.
\newblock Analysis of adaptive {BDF}2 scheme for diffusion equations.
\newblock {\em Math. Comp.}, 90(329):1207--1226, 2021.

\bibitem{LozinskiET:2009}
A.~Lozinski, M.~Picasso, and V.~Prachittham.
\newblock An anisotropic error estimator for the {C}rank--{N}icolson method:
  application to a parabolic problem.
\newblock {\em SIAM J. Sci. Comput.}, 31(4):2757--2783, 2009.

\bibitem{LubichET:2013}
C.~Lubich, D.~Mansour, and C.~Venkataraman.
\newblock Backward difference time discretization of parabolic differential
  equations on evolving surfaces.
\newblock {\em IMA J. Numer. Anal.}, 33(4):1365--1385, 2013.

\bibitem{MakridakisNochetto:2003}
C.~G. Makridakis and R.~H. Nochetto.
\newblock Elliptic reconstruction and a posteriori error estimates for
  parabolic problems.
\newblock {\em SIAM J. Numer. Anal.}, 41(4):1585--1594, 2003.

\bibitem{MakridakisNochetto:2006}
C.~G. Makridakis and R.~H. Nochetto.
\newblock A posteriori error analysis for higher order dissipative methods for
  evolution problems.
\newblock {\em Numer. Math.}, 104(4):489--514, 2006.

\bibitem{NevanlinnaOdeh:1981}
O.~Nevanlinna and F.~Odeh.
\newblock Multiplier techniques for linear multistep methods.
\newblock {\em Numer. Funct. Anal. Optim.}, 3(4):377--423, 1981.

\bibitem{Picasso:1998}
M.~Picasso.
\newblock Adaptive finite elements for a linear parabolic problem.
\newblock {\em Comput. Methods Appl. Mech. Engrg.}, 167(3--4):223--237, 1998.

\bibitem{ShampineReichelt:1997}
L.~F. Shampine and M.~W. Reichelt.
\newblock The {MATLAB} {ODE} suite.
\newblock {\em SIAM J. Sci. Comput.}, 18(1):1--22, 1997.

\bibitem{Smears:2025}
I.~Smears.
\newblock {A}n introduction to the a posteriori error analysis of parabolic
  partial differential equations.
\newblock In {\em Error control, adaptive discretizations, and applications,
  {P}art 4}, volume~61 of {\em Advances in Applied Mechanics}, pages 239--286.
  Elsevier, 2025.

\bibitem{Verfuerth:2003}
R.~Verf\"{u}rth.
\newblock A posteriori error estimates for finite element discretizations of
  the heat equation.
\newblock {\em Calcolo}, 40(3):195--212, 2003.

\bibitem{Verfuerth:2013}
R.~Verf\"{u}rth.
\newblock {\em A posteriori error estimation techniques for finite element
  methods}.
\newblock Numerical Mathematics and Scientific Computation. Oxford University
  Press, Oxford, 2013.

\bibitem{WangMaoHuang:2022}
W.~Wang, M.~Mao, and Y.~Huang.
\newblock Optimal a posteriori estimators for the variable step-size {BDF}2
  method for linear parabolic equations.
\newblock {\em J. Comput. Appl. Math.}, 413:114306, 2022.

\end{thebibliography}
\end{document}